\documentclass[11pt,a4paper,reqno]{amsart}
\usepackage{graphicx}
\usepackage{epsfig}
\usepackage{amsmath}
\usepackage{amsfonts}
\usepackage{amssymb}
\usepackage{amscd}

\numberwithin{equation}{section}

\begin{document}

\address{Misir J. Mardanov \newline
Ministry of Science and Education Republic of Azerbaijan,\newline
Institute of Mathematics and Mechanics,  Baku, Azerbaijan\newline
Baku State University, Baku , Azerbaijan}
\email{misirmardanov@yahoo.com}

\address{Telman K. Melikov \newline
Institute of Control Systems,\newline
Ministry of Science and Education Republic of Azerbaijan,\newline
Institute of Mathematics and Mechanics,  Baku, Azerbaijan\newline
}
\email{t.melik@rambler.ru }

\address{Gulnar V. Hajiyeva \newline
Ministry of Science and Education Republic of Azerbaijan,\newline
Institute of Mathematics and Mechanics,  Baku, Azerbaijan\newline
Baku Engineering University, Baku, Azerbaijan}
\email{gulnarhaciyeva607@gmail.com}

\begin{center}
{{\bf Necessary conditions for a minimum in variational problems with delay in the presence of degeneracies}}
\end{center}
\

\begin{center}
\textbf{M.J. Mardanov, T.K. Melikov, G.V. Hajiyeva }
\end{center}

\

\emph{This article  explores minimum of an extremal in the variational problem with delay under the degeneracy of the Weierstrass condition. Here for study the minimality of extremal, variations of the Weierstrass type are used in two forms: in the form of variations on the right with respect to the given point, and in the form of variations on the left with respect to the same point. Further, using these variations, formulas for the increments of the functional are obtained. The exploring of the minimality of the extremal with the help of these formulas is conducted under the assumption that the Weierstrass condition degenerates. As a result, considering different forms of degenerations (degeneracy of the Weieristrass condition at a single point and at points of a certain interval), we obtain the necessary conditions of the inequality type and the equality type for a strong and weak local minimum. A specific example is given to demonstrate the effectiveness of the results obtained in this article. }



\section{Introduction and problem statement}

In the present paper, we consider the following vector problem of variation with delayed argument:
\begin{equation} \label{GrindEQ__1_1_} 
S(x(\cdot ))=\int _{t_{0} }^{t_{1} }L(t,x(t),\, x(t-h),\, \dot{x}(t),\, \dot{x}(t-h))dt \to \mathop{\min }\limits_{x(\cdot )} ,                   
\end{equation} 
\begin{equation} \label{GrindEQ__1_2_} 
x(t)=\varphi(t),t\in [t_0-h,t_0],x(t_1)=x_1,\; x_1\in R^n,                         
\end{equation} 
here $R^{n} $ is $n$-dimensional space, $t_{0} ,\, t_{1} \in (-\infty ,+\infty ),\, \, x_{0} $ and $x_{1} $ are the given points $h={const}>0$,$t_1-t_0>h,$   functions  $L\left(t,x,y,\dot{x},\dot{y}\right):[t_0,t_1]\times R^n\times R^n\times R^n{\times R}^n\to R:=\left(-\infty,+\infty \right)$ and $\varphi (t)\in C^{1} \left([t_{0} -h,t_{0} ]\, ,\, R^{n} \right)$  are the given continuously differentiable functions with respect to   all their arguments and for each $(t,x,y,\dot{x},\dot{y})\in (t_1,+\infty )\times R^n\times R^n\times R^n\times R^n$  there is equality $L(t,x,y,\dot{x},\dot{y})=0$,where$\ \ \ \ y=y\left(t\right)=x\left(t-h\right),\ \ \ \dot{y}:=\dot{y}\left(t\right)=\dot{x}\left(t-h\right),\ \ t\in I:=[t_0,t_1]$, moreover $x(t)\in KC^1\left(\hat{I},R^n\right)$, where $\hat{I}\, :\, =\left[t_{0} -h,\, t_{1} \right]$   and $KC^{1} \left(\hat{I}\, ,\, R^{n} \right)$ is a class of piecewise-smooth function .

We call the functions $x(\cdot )\in KC^{1} \left(\hat{I}\, ,\, R^{n} \right)$ satisfying boundary conditions \eqref{GrindEQ__1_2_}, an admissible function.

Let us introduce the following notation to study the problem \eqref{GrindEQ__1_1_}, \eqref{GrindEQ__1_2_} below for compact and convenient notation along with the admissible function $\overline{x}(\cdot )$ :
\[\bar{L}(\tau ):=L(\tau ,\bar{x}(\tau ),\bar{y}(\tau ),\dot{\bar{x}}(\tau ),\dot{\bar{y}}(\tau ))\, , \] 
\begin{equation} \label{GrindEQ__1_3_} 
\bar{L}(\tau,\xi;\dot{\bar{x}}(\cdot )):=L(\tau,\bar{x}(\tau),\bar{y}(\tau),\dot{\bar{x}}(\tau)+\xi,\dot{\bar{y}}(\tau)), \bar{L}(\tau,\xi;\dot{\bar{y}}(\cdot )):=L(\tau,\bar{x}(\tau),\bar{y}(\tau),\dot{\bar{x}}(\tau),\dot{\bar{y}}(\tau)+\xi) , 
\end{equation} 
(the notations are defined similarly ${\bar{L}}_x(\tau),{\bar{L}}_y(\tau),{\bar{L}}_{\dot{x}}(\tau)$,$\bar{L}_{\dot{y}} (\tau )$ ${\bar{L}}_x(\tau,\xi;\dot{\bar{x}}(\cdot ))$ and    $\bar{L}_{y} (\tau ,\, \xi ;\, \dot{\bar{y}}(\cdot )\, )$) , where $\tau \in \{ t,t+h\} ,\, \, \, \tau \in I$,$\xi \in R^{n} $;
\begin{equation}\label{GrindEQ__1_4_}
E\left( \bar{L}\right) (\tau ,\xi ;\dot{\bar{x}}(\cdot )):=\bar{L}(\tau ,\xi
;\dot{\bar{x}}(t))-\bar{L}(\tau ){-\bar{L}}_{\dot{x}}^{T}(\tau )\xi ,
\end{equation}%
\[
E\left( \bar{L}\right) (\nu ,\xi ;\dot{\bar{y}}(\cdot )):=\bar{L}(\nu ,\xi ;%
\dot{\bar{y}}(\cdot ))-\bar{L}(\nu )-{\bar{L}}_{\dot{y}}^{T}(\nu)\xi ;
\]%
\[
Q_{k}\left( \bar{L}\right) (\tau ,\lambda ,\xi ;\dot{\bar{x}}(\cdot
)):=\lambda ^{k}E\left( \bar{L}\right) (\tau ,\xi ;\dot{\bar{x}}(\cdot
))+\left( 1-\lambda ^{k}\right) E\left( \bar{L}\right) \left( \tau ,\frac{%
\lambda }{\lambda -1}\xi ;\dot{\bar{x}}(\cdot )\right) ,
\]%
\begin{equation}\label{GrindEQ__1_5_}
Q_{k}\left( \bar{L}\right) (\nu ,\lambda ,\xi ;\dot{\bar{y}}(\cdot
)):=\lambda ^{k}E\left( \bar{L}\right) (\nu ,\xi ;\dot{\bar{y}}(\cdot
))+\left( 1-\lambda ^{k}\right) E\left( \bar{L}\right) \left( \nu ,\frac{%
\lambda }{\lambda -1}\xi ;\dot{\bar{y}}(\cdot )\right) ,
\end{equation}
where $\tau ,\nu \in \{t,t+h\}$, $\xi \in R^{n}$,$\lambda \in
(0,\,1),\,\,\,k=1,2$; 

\begin{eqnarray}\label{GrindEQ__1_6_}
M\left( {\bar{L}}_{x}\right) \left( \tau ,\lambda ,\xi ;\dot{\bar{x}}\left(
\cdot \right) \right)  &:&=\lambda \left[ {\bar{L}}_{x}^{T}\left( \tau ,\xi ;%
\dot{\bar{x}}\left( \cdot \right) \right) -{\bar{L}}_{x}^{T}\left( \tau\right) %
\right] \xi + \notag\\
&&+\left( 1-\lambda \right) \left[ {\bar{L}}_{x}^{T}\left( \tau ,\frac{%
\lambda }{\lambda -1}\xi ;\dot{\bar{x}}(\cdot )\right) -{\bar{L}}%
_{x}^{T}(\tau )\right] \xi , \\
M\left( {\bar{L}}_{y}\right) (\nu ,\lambda ,\xi ;\dot{\bar{y}}(\cdot ))
&:&=\lambda \left[ {\bar{L}}_{y}^{T}(\nu ,\xi ;\dot{\bar{y}}(\cdot ))-{\bar{L%
}}_{y}^{T}(\nu )\right] \xi +\notag \\
&&+\left( 1-\lambda \right) \left[ {\bar{L}}_{y}^{T}\ \left( \nu ,\frac{%
\lambda }{\lambda -1}\xi ;\dot{\bar{y}}(\cdot )\right) -{\bar{L}}%
_{y}^{T}(\nu )\right] \xi ,\notag
\end{eqnarray}
 where $\tau,\nu \in \{t,t+h\}$, $\xi \in R^{n}$,$\lambda
\in (0,1).$

Let us recall the concepts that were introduced, for example  in [22, 25]. The admissible function $\overline{x}(\cdot )$   is called a strong (weak) local minimum in the problem \eqref{GrindEQ__1_1_}, \eqref{GrindEQ__1_2_}, if there exists such a number $\bar{\delta }>0\, \, \, \left(\hat{\delta }>0\right)$  that the inequality $S(x(\cdot ))\ge S(\bar{x}(\cdot ))\, $ is fulfilled for all admissible functions $x(\cdot )$ for which
\[{\left\|x(\cdot )-\bar{x}(\cdot )\right\|}_{C\left(\hat{I},R^n\right)}\le \bar{\delta}{\left({max \left\{{\left\|x(\cdot )-\bar{x}(\cdot )\right\|}_{C\left(\hat{I},R^n\right)},{\left\|\dot{\bar{x}}(\cdot )-\dot{\bar{x}}(\cdot )\right\|}_{L_{\infty}\left(\hat{I},R^n\right)}\right\}\ }\le \hat{\delta}\right)}.\] 

In this case, we say that the admissible function $\bar{x}(\cdot )$ affords a strong (weak) local minimum in the problem \eqref{GrindEQ__1_1_}, \eqref{GrindEQ__1_2_} with $\bar{\delta }\, \, \, \left(\hat{\delta }\right)$ neighbourhood. Obviously, any strong local minimum at the same time is weak as well, but the opposite is not always true [23].

For the first time in 1970, Kamensky [12] investigated variational problems with delay and obtained an analogue of the Euler equation. It is easy to see that boundary condition \eqref{GrindEQ__1_2_} is more characteristic for the variational problem with delay compared to the boundary condition from [12]. Later, only in recent years the variational problem with delay have been studied in the works  [10,11,25]. 

Let us  also recall some information about the necessary conditions  for the minimality obtained in the works [10,25] . Namely:

(a)  if the admissible function $\bar{x}(\cdot )$ is a weak local minimum in the problem \eqref{GrindEQ__1_1_}, \eqref{GrindEQ__1_2_}, then along the function $\bar{x}(\cdot )$, the first variation of the problem \eqref{GrindEQ__1_1_}, \eqref{GrindEQ__1_2_} is equal to zero, i.e.  taking into account \eqref{GrindEQ__1_3_}, the following equality holds [10]:
\begin{equation} \label{GrindEQ__1_7_} 
\delta S(\delta x(\cdot );\bar{x}(\cdot ))=\int^{t_1}_{t_0}{\left\{\left[{\bar{L}}^T_x(t)+{\bar{L}}^T_y(t+h)\right]\delta x(t)+\left[{\bar{L}}^T_{\dot{x}}(t)+{\bar{L}}^T_{\dot{y}}(t+h)\right]\delta \dot{x}(t)\right\}dt}=0 ,  
\end{equation} 
      $\forall \delta x(t)\in KC^1\left(\hat{I},R^n\right)$, and $\delta x\left(t\right)=0,$ with $t\in [t_0-h,t_0]\cup \{t_1\}$;

(b) if the admissible function $\bar{x}(\cdot )$ is a weak local minimum in the problem \eqref{GrindEQ__1_1_}, \eqref{GrindEQ__1_2_} and the function $\ddot{\overline{x}}(\cdot )$ is continuous at the points in the set $\widetilde{{I}}\ \subset  I$ where $I \backslash \widetilde{I}-$ is a finite set, moreover, the functions $L\left(\cdot \right),\ \varphi\left(\cdot \right)$ are twice continuously differentiable with respect to the set of variables, then the function $\bar{x}(\cdot )$ is a solution of the Euler equation, i.e., the following equalities hold [10]:     
\begin{equation} \label{GrindEQ__1_8_} 
\frac{\ d}{dt}{\bar{L}}_{\dot{x}}\left(t\right)={\bar{L}}_x\left(t\right),\ \ \ t\in \widetilde{I}\cap (t_1-h,t_1],
\end{equation}                         
\[\frac{d}{dt}{[\bar{L}}_{\dot{x}}\left(t\right)+{\bar{L}}_{\dot{y}}\left(t+h\right)]={\bar{L}}_x\left(t\right)+{\bar{L}}_y\left(t+h\right),\ \ t\in \widetilde{{I}}\cap [t_0,t_1-h];\]

(c) if the admissible function $\bar{x}(\cdot )$ is strong local minimum in the problem \eqref{GrindEQ__1_1_}, \eqref{GrindEQ__1_2_}, then the analogue of the Weierstrass condition is fulfilled along it, i.e.,  for all $\xi \in R^{n} $ the following inequalities are valid [25]:
\begin{eqnarray}\label{GrindEQ__1_9_}
E\left( \bar{L}\right) (t,\xi ;\dot{\bar{x}}(\cdot )) &\ge&0,\forall t\in
I^{\ast }\cap[t_{1}-h,t_{1}], \notag\\
E\left( \bar{L}\right) (t,\xi ;\dot{\bar{x}}(\cdot ))+E\left( \bar{L}\right)
(t+h,\xi ;\dot{\bar{y}}(\cdot )) &\ge&0,\forall t\in I^{\ast }\cap[t_{0},t_{1}-h],
\end{eqnarray}
here $I^{*} \subseteq [t_{0} ,t_{1} ]$ is the set of points where the functions $\dot{\overline{x}}(\cdot )$ is continuous, in addition the functions $E\left(\bar{L}\right)\left(\cdot ;\dot{\bar{x}}\left(\cdot \right)\right)\ $ and  $E\left(\bar{L}\right)(\cdot ;\dot{\bar{y}}\left(\cdot \right))\ $ are defined by \eqref{GrindEQ__1_4_};

  (d) if the admissible function $\bar{x}(\cdot )$ is weak local minimum in the problem \eqref{GrindEQ__1_1_} \eqref{GrindEQ__1_2_}, then there exists a number $\delta >0$ such that for every $\xi \in B_\delta(0)$ the inequality \eqref{GrindEQ__1_9_} holds, where symbol $B_\delta(0)$ is a ball of the radius $\delta $ centered at the point $0\in R^n$  [25].

\textbf{Remark 1.1.} Following [3, 23], the minimum conditions \eqref{GrindEQ__1_9_} are valid on the right and on the left at discontinuity points of the function $\dot{\bar{x}}(\cdot )$.

We call a solution of the problem \eqref{GrindEQ__1_1_}, \eqref{GrindEQ__1_2_} as an extremal in the problem \eqref{GrindEQ__1_1_}, \eqref{GrindEQ__1_2_}, if it satisfies equation \eqref{GrindEQ__1_7_}.

It should be noted that the application of the Weierstrass condition \eqref{GrindEQ__1_9_} as a necessary condition for minimum is more effective if at each point $t\in I$  , the inequality \eqref{GrindEQ__1_9_} turns into equality only at a single point $\xi=0$. However, it may happen that at least at one point $\theta\in I$ inequality \eqref{GrindEQ__1_9_} turns into equality at several points $\xi\in R^n$. In this case, as a rule, it is usually said that the Weierstrass conditions \eqref{GrindEQ__1_9_} degenerates at the point $\theta $.

It is evident that the problem \eqref{GrindEQ__1_1_}, \eqref{GrindEQ__1_2_} in the terms of optimal control theory takes the following form:
\begin{equation} \label{GrindEQ__1_10_} 
S\left(x\left(\cdot \right),x_{n+1}\left(\cdot \right)\right)=x_{n+1}\left(t_1\right)\to {\mathop{\min }_{u(\cdot )} ,\ } 
\end{equation} 
\begin{equation} \label{GrindEQ__1_11_} 
\left\{ \begin{array}{c}
\dot{x}\left(t\right)=u\left(t\right),\ \ u\left(t\right)\in KC\left(\hat{{I}},R^n\right)\ , \ \ \ \ \ \ \ \ \ \ \ \ \ \ \ \ \ \ \ \ \ \ \ \ \  \\ 
{\dot{x}}_{n+1}\left(t\right)=L\left(t,x\left(t\right),x\left(t-h\right),u\left(t\right),u\left(t-h\right)\right),t\in\left(t_0,t_1\right], \end{array}
\right. 
\end{equation} 
\begin{equation} \label{GrindEQ__1_12_} 
x\left(t\right)=\varphi\left(t\right),u\left(t\right)=\dot{\varphi}\left(t\right),t\in\left[t_0-h,t_0\right],x\left(t_1\right)=x_1,x_{n+1}\left(t_0\right)=0, 
\end{equation}
where  $KC\left(\hat{{ I}},R^n\right)$ - the class of piecewise continuous functions.

As is known, degenerate cases in optimal control theory are studied in terms of singular controls. For the optimality of singular controls in control problems described by ordinary differential equations, several important results have been obtained (see, for example, [1,4,5,7,13-15,24,29,30]). Subsequently, some of these results have been significantly generalized for control problems with delays only in the forms, namely, for the problem with delays in phase variables (see, for example, [2,18, 26-29]) and for problems with delays in controls (see for example [16,17,20]).

It should be noted that, in the problems of the form \eqref{GrindEQ__1_10_} - \eqref{GrindEQ__1_12_} the optimality of singular controls has not been studied, since in the problem \eqref{GrindEQ__1_10_} - \eqref{GrindEQ__1_12_} the phase variable and also the control variable depend on the delay. Therefore, considering the above, it can be argued  that the study of the problem \eqref{GrindEQ__1_10_}-\eqref{GrindEQ__1_12_}  i.e. \eqref{GrindEQ__1_1_}, \eqref{GrindEQ__1_2_} in the presence of degeneracy is actual.

The main aim of our paper  is to minimize the extremal in the problem \eqref{GrindEQ__1_1_}, \eqref{GrindEQ__1_2_} under the degeneracy of the Weierstrass condition \eqref{GrindEQ__1_9_}. Developing the method developed in [23], we obtain the necessary equality-type and inequality-type conditions for a strong as well as a weak local minimum. The derivation of these minimum conditions is a significant generalization of the corresponding statements in [23].

The structure of the present paper is outlined by the following scheme. In the second section of the article, special variations are introduced for the extremal of problem \eqref{GrindEQ__1_1_}, \eqref{GrindEQ__1_2_} (see (2.1) and \eqref{GrindEQ__2_3_}), and lemmas are proved (see Lemmas (2.1) and \eqref{GrindEQ__2_2_}). In the third and fourth sections, using variation \eqref{GrindEQ__2_1_} and \eqref{GrindEQ__2_3_}, formulas for the increment of the functional \eqref{GrindEQ__1_1_} are obtained, namely, Lemmas 3.1, 3.2, 4.1 and 4.2 are proved. In the fifth and sixth sections, considering the degeneration of the Weierstrass condition \eqref{GrindEQ__1_9_}, based on the increment formulas obtained in the sections 3 and 4, theorems are proved (see Theorems 5.2, 6.1 and 6.2) on the necessary conditions for a minimum. Finally, in the seventh section, a specific example is given that shows the effectiveness of the statements of Theorem 5.1 and the prospects for further development and generalization of this work are noted.

\section{Special variations of the extremal of the problem (1.1),(1.2)  and the lemma.}

Let an admissible function $\bar{x}(\cdot )$ be an extremal of the problem \eqref{GrindEQ__1_1_}, \eqref{GrindEQ__1_2_} and $\vartheta:=\left(\theta,\lambda,\xi\right)\in \left[t_0,t_1-h\right]\times (0,1)\times R^n\backslash \{0\}$  -be an arbitrary fixed point.

Following [23], we take into consideration the following special variations of the extremal $\bar{x}(\cdot )$:
\begin{enumerate}
\item  variation introduced on the right with respect to the point $\theta \in [t_{0} ,t_{1} -h)$  :
\begin{equation} \label{GrindEQ__2_1_} 
x^{(+)} (t;\vartheta ,\varepsilon )=\bar{x}(t)+q^{(+)} (t;\vartheta ,\varepsilon ),\, \, \, \, \, t\in \hat{I}, 
\end{equation} 
\end{enumerate}
where the function $q^{(+)} (t;\vartheta ,\varepsilon )$is defined as:
\begin{equation} \label{GrindEQ__2_2_} 
q^{(+)}(t;\vartheta,\varepsilon)=\left\{ \begin{array}{c}
(t-\theta)\xi,\ \ \ \ \ \ \ \ \ \ \ \ \ \ \ \ \ t\in [\theta,\theta+\lambda \varepsilon), \\ 
\frac{\lambda}{\lambda-1}(t-\theta-\varepsilon)\xi,\ \ \ t\in [\theta+\lambda \varepsilon,\theta+\varepsilon), \\ 
0,\ \ \ \ \ \ \ \ \ \ \ \ \ \ \ \ \ \ \ \ \ \ \ \ \ \ \ \ \ \ \ \ t\in \hat{I}\backslash [\theta,\theta+\varepsilon), \end{array}
\right. 
\end{equation} 
here $\lambda \in (0,\, 1),\, \, \, \, \xi \in R^{n} \backslash \{ 0\} $, $\vartheta =(\theta ,\lambda ,\xi )$ and  $\varepsilon\in (0,\overline{\varepsilon})$,  $\overline{\varepsilon}={min \{\ }h,t_1-\theta-h\}$;

(ii) variation introduced on the left with respect to the point $\theta \in (t_{0} ,t_{1} -h]$ :
\begin{equation} \label{GrindEQ__2_3_} 
x^{(-)} (t;\vartheta ,\varepsilon )=\bar{x}(t)+q^{(-)} (t;\vartheta ,\varepsilon ),\, \, \, \, \, t\in \hat{I},   
\end{equation} 
where the function $q^{(-)} (t;\vartheta ,\varepsilon )$ is defined as:
\begin{equation}\label{GrindEQ__2_4_}
q^{(-)}(t;\vartheta ,\varepsilon )=\left\{ 
\begin{array}{c}
(t-\theta )\xi ,\ \ \ \ \ \ \ \ \ \ \ \ \ \ \ \ t\in (\theta -\lambda
\varepsilon ,\theta ], \\ 
\frac{\lambda }{\lambda -1}(t-\theta +\varepsilon )\xi ,\; t\in (\theta
-\varepsilon ,\theta -\lambda \varepsilon ], \\ 
0,\ \ \ \ \ \ \ \ \ \ \ \ \ \ \ \ \ \ \ \ \ \ \ \ \ \ \ \ \ t\in \hat{I}%
\backslash (\theta -\varepsilon ,\theta ],%
\end{array}%
\right.   
\end{equation}
here $\lambda \in (0,\, 1),\, \, \, \, \xi \in R^{n} \backslash \{ 0\} $, $\vartheta =(\theta ,\lambda ,\xi )$ and  $\tilde{\varepsilon }=\min \{ h,\theta -t_{0} \} $.

It is obvious that for any $\varepsilon \in (0,\hat{\varepsilon })$,  where $\hat{\varepsilon}={min \{\ }\bar{\varepsilon},\tilde{\varepsilon}\}$, functions  $q^{(+)} (\cdot \, \, ;\vartheta ,\varepsilon )$ and  $q^{(-)}(\cdot \vartheta,\varepsilon)$ are elements of the space $KC^1\left(\hat{I},R^n\right)$  and their derivatives $\dot{q}^{(+)} (\, \cdot \, \, ;\vartheta ,\varepsilon )$ and ${\dot{q}}^{(-)}(\cdot ;\vartheta ,\varepsilon)$ are calculated using the following  formulas:

   \begin{equation}\label{GrindEQ__2_5_}
{\dot{q}}^{(+)}(t;\vartheta ,\varepsilon )=\left\{ 
\begin{array}{c}
\xi,\ \ \ \ \ \ \ \ t\in \lbrack \theta ,\theta +\lambda
\varepsilon ], \\ 
\frac{\lambda }{\lambda -1}\xi ,\ t\in \lbrack \theta +\lambda \varepsilon
,\theta +\varepsilon \ ], \\ 
0,\ \ \ \ \ \ \ \ t\in \hat{I}\backslash (\theta ,\theta +\varepsilon ),%
\end{array}%
\right.   
\end{equation}%
\begin{equation}\label{GrindEQ__2_6_}
{\dot{q}}^{(-)}(t;\vartheta ,\varepsilon )=\left\{ 
\begin{array}{c}
\xi ,\ \ \ \ \ \ \ t\in \lbrack \theta -\lambda \varepsilon ,\theta ], \\ 
\frac{\lambda }{\lambda -1}\xi ,\ t\in \lbrack \theta -\xi ,\theta -\lambda
\varepsilon ], \\ 
0,\ \ \ \ \ \ \ \ t\in \hat{I}\backslash (\theta -\varepsilon, \theta ).
\end{array}%
\right.   
\end{equation}

As seen, the derivative $\dot{q}^{(+)} (\, \cdot ;\vartheta ,\varepsilon )$ is calculated both on the right and on the left at the points  $\theta,\theta+\lambda \varepsilon$ and $\theta+\varepsilon$ and the derivative $\dot{q}^{(-)} (t;\vartheta ,\varepsilon )$ is calculated on both the right and on the left at the points $\theta-\varepsilon, \theta-\lambda \varepsilon,$ and $\theta$.

\textbf{Lemma 2.1.} Let an admissible function $\bar{x}(\cdot )$ be an extremal of the problem \eqref{GrindEQ__1_1_}, \eqref{GrindEQ__1_2_}. Then for each $\vartheta =(\theta ,\lambda ,\xi )\in [t_{0} ,t_{1} -h)\times (0,\, 1)\times R^{n} \backslash \{ 0\} $ and for all   $\varepsilon \in (0,\bar{\varepsilon })$ the equality holds: 
\[
\delta S(q^{(+)}(\cdot ;\vartheta ,\varepsilon );\bar{x}(\cdot
))=\int_{\theta }^{\theta +\lambda \varepsilon }{\left[ \left( {\bar{L}}%
_{x}^{T}(t)+{\bar{L}}_{y}^{T}(t+h)\right) (t-\theta )+{\bar{L}}_{\dot{x}%
}^{T}(t)+{\bar{L}}_{\dot{y}}^{T}(t+h)\right] \xi dt}+
\]%
\begin{equation}\label{GrindEQ__2_7_}
+\frac{\lambda }{\lambda -1}\int_{\theta +\lambda \varepsilon }^{\theta
+\varepsilon }{\left[ \left( {\bar{L}}_{x}^{T}\left( t\right) +{\bar{L}}%
_{y}^{T}\left( t+h\right) \right) (t-\theta -\varepsilon )+{\bar{L}}_{\dot{x}%
}^{T}(t)+{\bar{L}}_{\dot{y}}^{T}(t+h)\right] \xi dt}=0,
\end{equation} 
where $q^{(+)} (\, \cdot \, \, ;\vartheta ,\varepsilon )$ is determined by \eqref{GrindEQ__2_2_}.

Proof of  Lemma 2.1 follows \eqref{GrindEQ__1_7_} taking into account $\delta \, x(\, \cdot \, )=q^{(+)} (\, \cdot \, \, ;\vartheta ,\varepsilon )$ and $\delta \, \dot{x}(\, \cdot \, )=\dot{q}^{(+)} (\, \cdot \, \, ;\vartheta ,\varepsilon )$, where $q^{(+)} (\, \cdot \, \, ;\vartheta ,\varepsilon )$ and $\dot{q}^{(+)} (\, \cdot \, \, ;\vartheta ,\varepsilon )$ are determined by \eqref{GrindEQ__2_2_} and \eqref{GrindEQ__2_5_}, respectively.

\textbf{Lemma 2.2.} Let an admissible function $\bar{x}(\cdot )$ be an extremal of the problem \eqref{GrindEQ__1_1_}, \eqref{GrindEQ__1_2_}. Then for each $\vartheta =(\theta ,\lambda ,\xi )\in (t_{0} ,t_{1} -h]\times (0,\, 1)\times R^{n} \backslash \{ 0\} $ and for all $\varepsilon \in (0,\tilde{\varepsilon })$ the equality holds:
\[
\delta S(q^{(-)}(\cdot ;\vartheta ,\varepsilon );\bar{x}(\cdot
))=\int_{\theta -\lambda \varepsilon }^{\theta }{\left[ \left( {\bar{L}}%
_{x}^{T}(t)+{\bar{L}}_{y}^{T}(t+h)\right) (t-\theta )+{\bar{L}}_{\dot{x}%
}^{T}(t)+{\bar{L}}_{\dot{y}}^{T}(t+h)\right] \xi dt}+
\]%
\begin{equation} \label{GrindEQ__2_8_}
+\frac{\lambda }{\lambda -1}\int_{\theta -\varepsilon }^{\theta -\lambda
\varepsilon }{\left[ \left( {\bar{L}}_{x}^{T}(t)+{\bar{L}}%
_{y}^{T}(t+h)\right) (t-\theta +\varepsilon )+{\bar{L}}_{\dot{x}}^{T}(t)+{
\bar{L}}_{\dot{y}}^{T}(t+h)\right] \xi dt}=0, 
\end{equation}
where $q^{(-)} (\, \cdot \, \, ;\vartheta ,\varepsilon )$ is determined by \eqref{GrindEQ__2_4_}.

Proof of Lemma 2.2 follows \eqref{GrindEQ__1_7_} taking into account\textbf{ }$\delta \, x(\, \cdot \, )=q^{(-)} (\, \cdot \, \, ;\vartheta ,\varepsilon )$ and $\delta \, \dot{x}(\, \cdot \, )=\dot{q}^{(-)} (\, \cdot \, \, ;\vartheta ,\varepsilon )$, where $q^{(-)} (\, \cdot \, \, ;\vartheta ,\varepsilon )$ and $\dot{q}^{(-)} (\cdot \,  ;\vartheta ,\varepsilon )$ are determined by \eqref{GrindEQ__2_4_} and \eqref{GrindEQ__2_6_}, respectively.

\section{ Formula for the increment of the functional in the problem (1.1),(1.2)  on variations (2.1)} 

Let  $\bar{x}(\cdot )$ be an extremal of the problem \eqref{GrindEQ__1_1_}, \eqref{GrindEQ__1_2_} and suppose that the derivative function of  $\overline{x}\left(\cdot \right)\ \ $is continuous on the right at the points $\theta -h,\, \, \theta ,\, \theta +h$. Continuing of the study, let's calculate the increment:
\begin{equation} \label{GrindEQ__3_1_} 
S(x^{(+)} (\, \cdot \, ;\, \vartheta ,\varepsilon ))-S(\bar{x}(\cdot ))\, \, =\, :\Delta \, S^{(+)} (q^{(+)} (\, \cdot \, ;\, \vartheta ,\varepsilon );\bar{x}(\cdot )) 
\end{equation} 
functional \eqref{GrindEQ__1_1_} in the problem \eqref{GrindEQ__1_1_}, \eqref{GrindEQ__1_2_} on variations \eqref{GrindEQ__2_1_} up to $o(\varepsilon ^{2} )$, where $\vartheta=(\theta,\lambda,\xi),\theta \in [t_0,t_1-h)$, $\varepsilon \in (0,\bar{\varepsilon })$, the number $\bar{\varepsilon }>0$ is defined in  \eqref{GrindEQ__2_2_}, $o(\varepsilon^2)/\varepsilon^2\to 0$,  at  $\varepsilon\to+0$. 

Increment \eqref{GrindEQ__3_1_}, considering \eqref{GrindEQ__2_1_}, \eqref{GrindEQ__2_2_}, \eqref{GrindEQ__2_5_} and Lemma 2.1 takes the following form:
\[
\Delta ^{(+)}S(q^{(+)}(\cdot \,;\,\vartheta ,\varepsilon );\bar{x}(\cdot
))\,\,=\,\Delta ^{(+)}\,S_{1}(q^{(+)}(\cdot \,;\,\vartheta ,\varepsilon );%
\bar{x}(\cdot ))+ 
\]%
\begin{equation}\label{GrindEQ__3_2_}
+\Delta ^{(+)}S_{2}(q^{(+)}(\cdot ;\vartheta ,\varepsilon );\bar{x}(\cdot
))-\delta S(q^{(+)}(\cdot ;\vartheta ,\varepsilon);\bar{x}(\cdot )),
\end{equation}%
here 
\[
\Delta ^{(+)}S_{1}(q^{(+)}(\,\cdot \,;\,\vartheta ,\varepsilon );\bar{x}%
(\cdot ))\,\,=\int_{\theta }^{\theta +\varepsilon }\,\left[ L\left( t,\bar{x}%
(t),\bar{y}(t),\dot{\bar{x}}(t)+\dot{q}_{\varepsilon }^{(+)}(t),\,\dot{\bar{y%
}}(t)\right) -\bar{L}(t)\right] \,dt\,+\, 
\]%
\begin{equation}\label{GrindEQ__3_3_}
+\int_{\theta +h}^{\theta +h+\varepsilon }\,\left[ L\left( t,\bar{x}(t),\bar{%
y}(t),\dot{\bar{x}}(t),\,\dot{\bar{y}}(t)+\dot{q}_{\varepsilon
}^{(+)}(t-h)\right) -\bar{L}(t)\right] \,dt,  
\end{equation}%
\[
\Delta ^{(+)}S_{2}(q^{(+)}(\cdot \,;\,\vartheta ,\varepsilon );\bar{x}(\cdot
))=\int_{\theta }^{\theta +\varepsilon}{\left[ L\left( t,\bar{x}(t)+q_{\varepsilon
}^{(+)}(t),\bar{y}(t),\dot{\bar{x}}(t)+{\dot{q}}_{\varepsilon }^{(+)}(t),%
\dot{\bar{y}}(t)\right) \right. }- 
\]%
\[
\left. -L\left( t,\bar{x}(t),\bar{y}(t),\dot{\bar{x}}(t)+\dot{q}%
_{\varepsilon }^{(+)}(t),\,\dot{\bar{y}}(t)\right) \,\right] \,dt+ 
\]%
\[
+\int_{\theta +h}^{\theta +h+\varepsilon }{\left[ L\left( t,\bar{x}(t),\bar{y%
}(t)+q_{\varepsilon }^{(+)}(t-h),\dot{\bar{x}}(t),\dot{\bar{y}}(t)+{\dot{q}}%
_{\varepsilon }^{(+)}(t-h)\right) \right. }- 
\]%
\begin{equation}\label{GrindEQ__3_4_}
\left. -L\left( t,\bar{x}(t),\bar{y}(t),\dot{\bar{x}}(t),\,\dot{\bar{y}}(t)+%
\dot{q}_{\varepsilon }^{(+)}(t-h)\right) \,\right] \,dt,
\end{equation}%
where $q_{\varepsilon}^{(+)}(\cdot ):=q^{(+)}(\cdot ;\,\vartheta ,\varepsilon )$, ${%
\dot{q}}_{\varepsilon}^{(+)}(\cdot ):={\dot{q}}^{(+)}(\cdot ;\,\vartheta ,\varepsilon
),\varepsilon \in (0,\bar{\varepsilon })$, then $\delta S(q^{\left( +\right)
}\left( \cdot ;\,\vartheta ,\varepsilon \right) ;\bar{x}\left( \cdot \right)
)$ is defined by \eqref{GrindEQ__2_7_}.

From \eqref{GrindEQ__3_3_} and \eqref{GrindEQ__3_4_}, considering %
\eqref{GrindEQ__2_2_}, \eqref{GrindEQ__2_5_} and notation %
\eqref{GrindEQ__1_3_}, we obtain: 
\[
\Delta ^{(+)}S_{1}(q^{(+)}(\cdot ;\,\vartheta ,\varepsilon );\bar{x}(\cdot
))=\int_{\theta }^{\theta +\lambda \varepsilon }{\left\{ \left[ \bar{L}%
(t,\xi ;\dot{\bar{x}}(\cdot ))-\bar{L}(t)\right] +{\left. \left[ \bar{L}%
(\tau ,\xi ;\dot{\bar{y}}(\cdot ))-\bar{L}(\tau)\right] \right\vert }_{\tau
=t+h}\right\} dt}+
\]%
\begin{equation}\label{GrindEQ__3_5_}
+\int_{\theta +\lambda \varepsilon }^{\theta +\varepsilon }{\left\{ \left[ 
\bar{L}\left( t,\frac{\lambda }{\lambda -1}\xi ;\dot{\bar{x}}(\cdot )\right)
-\bar{L}(t)\right] +{\left. \left[ \bar{L}\left( \tau ,\frac{\lambda }{%
\lambda -1}\xi ;\dot{\bar{y}}(\cdot )\right) -\bar{L}(\tau)\right] \right\vert }%
_{\tau =t+h}\right\} dt},  
\end{equation}%
\[
\Delta ^{(+)}S_{2}(q^{(+)}(\,\cdot \,;\,\vartheta ,\varepsilon );\bar{x}%
(\cdot ))\,\,=\int_{\theta }^{\theta +\varepsilon }\,\,\left\{ \left[
L\left( t,\bar{x}(t)+(t-\theta )\xi ,\bar{y}(t),\dot{\bar{x}}(t)+\xi ,\,\dot{%
\bar{y}}(t)\right) \right. \right. -\,
\]%
\[
\left. -L\left( t,\bar{x}(t),\bar{y}(t),\dot{\bar{x}}(t)+\xi ,\dot{\bar{y}}%
(t)\right) \right] +\left[ L\left( \tau ,\bar{x}(\tau ),\bar{y}(\tau
)+(t-\theta )\xi ,\dot{\bar{x}}(\tau ),\dot{\bar{y}}(\tau )+\xi \right)
\right. -
\]%
\[
{\left. {\left. \left. -L\left( \tau ,\bar{x}(\tau ),\bar{y}(\tau ),\dot{%
\bar{x}}(\tau ),\dot{\bar{y}}(\tau )+\xi \right) \right] \right\vert }_{\tau
=t+h}\right\} }dt+
\]%
\[
+\int_{\theta +\lambda \varepsilon }^{\theta +\varepsilon }\,\left\{ \,\left[
L\left( t,\bar{x}(t)+\frac{\lambda }{\lambda -1}(t-\theta -\varepsilon )\xi ,%
\bar{y}(t),\dot{\bar{x}}(t)+\frac{\lambda }{\lambda -1}\xi ,\dot{\bar{y}}%
(t)\right) -\right. \right. 
\]%
\[
\left. -L\left( t,\bar{x}(t),\bar{y}(t),\dot{\bar{x}}(t)+\frac{\lambda }{%
\lambda -1}\xi ,\dot{\bar{y}}(t)\right) \,\right] +
\]%
\[
+\left[ L\left( \tau ,\bar{x}(\tau ),\bar{y}(\tau )+\frac{\lambda }{\lambda
-1}(t-\theta -\varepsilon )\xi ,\dot{\bar{x}}(\tau),\dot{\bar{y}}(\tau)+\frac{%
\lambda }{\lambda -1}\xi \right) -\right. 
\]
\begin{equation}\label{GrindEQ__3_6_}
\left. {\left. \left. -L\left( \tau,\bar{x}(\tau),\bar{y}(\tau),\dot{\bar{x}}(\tau),\dot{%
\bar{y}}(\tau)+\frac{\lambda }{\lambda -1}\xi \right) \right] \right\vert }%
_{\tau=t+h}\right\} dt.  
\end{equation}

According to the Taylor formula, formula \eqref{GrindEQ__3_6_} takes the following form:
\[
\Delta ^{\left( +\right) }S_{2}(q^{\left( +\right) }(\cdot ;\,\vartheta
,\varepsilon );\bar{x}\left( \cdot \right) )=\int_{\theta }^{\theta +\lambda
\varepsilon }{\left\{ \left[ {\bar{L}}_{x}^{T}(t,\xi ;\dot{\bar{x}}(\cdot ))+%
{\bar{L}}_{y}^{T}(t+h,\xi ;\dot{\bar{y}}(\cdot ))\right] (t-\theta )\xi
+o(t-\theta)\right\} dt}+
\]%
\[
+\int_{\theta +\lambda \varepsilon }^{\theta +\varepsilon }{\left\{ \left[ {%
\bar{L}}_{x}^{T}\left( t,\frac{\lambda }{\lambda -1}\xi ;\dot{\bar{x}}(\cdot
)\right) +{\bar{L}}_{y}^{T}\left( t+h,\frac{\lambda }{\lambda -1}\xi ;\dot{%
\bar{y}}(\cdot )\right) \right] \xi (t-\theta-\varepsilon )+o(t-\theta
-\varepsilon )\right\} dt=}
\]%
\[
=\frac{1}{2}\lambda ^{2}\varepsilon ^{2}\left[ {\bar{L}}_{x}^{T}(\theta
_{+},\xi ;\dot{\bar{x}}(\cdot ))+{\bar{L}}_{y}^{T}({(\theta +h)}_{+},\xi;%
\dot{\bar{y}}(\cdot ))\right] \xi +
\]%
\begin{equation}\label{GrindEQ__3_7_}
+\frac{\varepsilon ^{2}\lambda (1-\lambda )}{2}\left[ {\bar{L}}%
_{x}^{T}\left( \theta _{+},\frac{\lambda }{\lambda -1}\xi ;\dot{\bar{x}}%
(\cdot )\right) +{\bar{L}}_{y}^{T}\left( {(\theta +h)}_{+},\frac{\lambda }{%
\lambda -1}\xi ;\dot{\bar{y}}(\cdot )\right) \right] \xi +o(\varepsilon
^{2}).  
\end{equation}

In addition, from \eqref{GrindEQ__2_7_} we have 
\[
\delta S(q^{(+)}(\cdot ;\,\vartheta ,\varepsilon );\bar{x}\left( \cdot \right)
)=\int_{\theta }^{\theta +\lambda \varepsilon }{\left[ {\bar{L}}_{\dot{x}%
}^{T}(t)+{\bar{L}}_{\dot{y}}^{T}(t+h)\right] \xi dt}+
\]%
\[
+\frac{\lambda }{\lambda -1}\int_{\theta +\lambda \varepsilon }^{\theta
+\varepsilon }{\left[ {\bar{L}}_{\dot{x}}^{T}\left( t\right) +{\bar{L}}_{%
\dot{y}}^{T}\left( t+h\right) \right] \xi dt+\frac{1}{2}\lambda
^{2}\varepsilon ^{2}}\left[ {\bar{L}}_{x}^{T}\left( \theta _{+}\right) +{%
\bar{L}}_{y}^{T}\left( {(\theta +h)}_{+}\right) \right] \xi +
\]%
\begin{equation}\label{GrindEQ__3_8_}
+\frac{1}{2}\varepsilon ^{2}\lambda (1-\lambda )\left[ {\bar{L}}%
_{x}^{T}\left( \theta _{+}\right) +{\bar{L}}_{y}^{T}\left( {(\theta +h)}%
_{+}\right) \right] \xi +o(\varepsilon ^{2}).  
\end{equation}

When deriving \eqref{GrindEQ__3_7_} and \eqref{GrindEQ__3_8_}, we considered  into the assumptions about the right continuity of the function $\dot{\overline{x}}\left(\cdot \right)\ $ at the points $\theta-h,\theta,\theta+h$.

We consider \eqref{GrindEQ__3_5_}, \eqref{GrindEQ__3_7_} and \eqref{GrindEQ__3_8_} in \eqref{GrindEQ__3_2_}. As a result, considering notations \eqref{GrindEQ__1_3_}, \eqref{GrindEQ__1_4_} and \eqref{GrindEQ__1_6_} we obtain:
\[
\Delta ^{(+)}S(q^{(+)}(\cdot ;\,\vartheta ,\varepsilon );\overline{x}\
):=\int_{\theta }^{\theta +\lambda \varepsilon }{\left[ E\left( \bar{L}%
\right) (t,\xi ;\dot{\bar{x}}(\cdot ))+E\left( \bar{L}\right) (t+h,\xi ;\dot{%
\bar{y}}(\cdot ))\right] dt}+
\]%
\[
+\int_{\theta +\lambda \varepsilon }^{\theta +\varepsilon }{\left[ E\left( 
\bar{L}\right) \left( t,\frac{\lambda }{\lambda -1}\xi ;\dot{\bar{x}}(\cdot
)\right) +E\left( \bar{L}\right) \left( t+h,\frac{\lambda }{\lambda -1}\xi ;
\dot{\bar{y}}(\cdot )\right) \right] dt+}
\]%
\begin{equation}\label{GrindEQ__3_9_}
 +\frac{1}{2}\xi ^{2}\lambda \left[ M\left( {\bar{%
L}}_{x}\right) \left( \theta _{+},\lambda ,\xi ;\dot{\bar{x}}(\cdot )\right)
+M\left( {\bar{L}}_{y}\right) \left( {(\theta +h)}_{+},\lambda ,\xi ;\dot{%
\bar{y}}(\cdot )\right) \right] +o(\varepsilon^{2}),  
\end{equation}%
where functions $E\left( \bar{L}\right) \left( \cdot \,\,,\,\cdot \,\,;\,%
\dot{\bar{x}}(\cdot )\right) $ and $E\left( \bar{L}\right) \left( \cdot
\,\,,\,\cdot \,\,;\,\dot{\bar{y}}(\cdot )\right) $ are determined by %
\eqref{GrindEQ__1_4_}, and functions $M\left( \bar{L}\right) \left( \cdot
\,\,,\,\cdot \,\,;\,\dot{\bar{x}}(\cdot )\right) $ and $M\left( \bar{L}%
\right) \left( \cdot \,\,,\,\cdot \,\,;\,\dot{\bar{y}}(\cdot )\right) $ are
determined by \eqref{GrindEQ__1_6_}.

Thus, the following statement is proved.

\textbf{Lemma 3.1.} Let the functions $L(\, \cdot \, )$ and $f\left(\cdot \right)$ be continuously differentiable with respect to the totality of their arguments and the admissible function $\bar{x}(\cdot )$ is a solution of equation \eqref{GrindEQ__1_7_}. In addition, let the function  $\dot{\overline{x}}\left(\cdot \right)  $ be continuous on the right at the points $\theta-h,\ \theta $ and $\theta+h$, where $\theta \in [t_0,t_1-h)$.  Then, for sufficiently small values $\varepsilon >0$, the expansion \eqref{GrindEQ__3_9_} holds.

Based on \eqref{GrindEQ__3_9_}, the following Lemma is proven.

\textbf{Lemma 3.2.} Let functions $L(\cdot )$, $L_{\dot{x}} (\cdot )$ , $L_{\dot{y}} (\cdot )$ be continuous differentiable with respect to the totality of their arguments and the function $\varphi\left(\cdot \right)-$ twice continuously differentiable. In addition, let an admissible function $\bar{x}(\cdot )$ is an extremal of the problem \eqref{GrindEQ__1_1_}, \eqref{GrindEQ__1_2_} and let it be twice continuously differentiable on the right in the  semi-neighborhoods of the points $\theta -h,\, \theta $ and  $\theta +h$. Then for sufficiently small values $\varepsilon >0$ there is a representation of the form:
\[
\Delta ^{(+)}S(q^{(+)}(\cdot ;\,\vartheta ,\varepsilon );\overline{x}(\cdot
))=\varepsilon \left[ Q_{1}\left( \bar{L}\right) (\theta _{+},\lambda ,\xi ;%
\dot{\bar{x}}(\cdot ))+Q_{1}((\theta +h)_{+},\lambda ,\xi ;\dot{\bar{y}}%
(\cdot ))\right] +
\]%
\[
+\frac{1}{2}\varepsilon ^{2}\left\{ \lambda \left[ M\left( {\bar{L}}%
_{x}\right) \left( \theta _{+},\lambda ,\xi ;\dot{\bar{x}}(\cdot )\right)
+M\left( {\bar{L}}_{y}\right) \left( (\theta +h)_{+},\lambda ,\xi ;\dot{\bar{%
y}}(\cdot )\right) \right] \right. +
\]%
\begin{equation}\label{GrindEQ__3_10_}
\left. +\frac{d}{dt}\left[ Q_{2}\left( \bar{L}\right) (\theta _{+},\lambda
,\xi ;\dot{\bar{x}}(\cdot ))+Q_{2}((\theta +h)_{+},\lambda ,\xi ;\dot{\bar{y}%
}(\cdot ))\right] \right\} +o(\varepsilon ^{2}),  
\end{equation}
where $Q_{k} \left(\bar{L}\right)\left(\, \cdot \, ;\, \dot{\bar{x}}(\cdot )\right)$ and $Q_{k} \left(\bar{L}\right)\left(\, \cdot \, ;\, \dot{\bar{y}}(\cdot )\right),\, k=1,2$ are defined by \eqref{GrindEQ__1_5_} and $M\left(\bar{L}_{x}^{} \right)\left(\, \cdot \, ;\, \dot{\bar{x}}(\cdot )\right)$ and $M\left(\bar{L}_{y}^{} \right)\left(\, \cdot \, ;\, \dot{\bar{y}}(\cdot )\right)$  by \eqref{GrindEQ__1_6_}.

\textbf{Proof.} Using the Taylor's formula, considering the smoothness assumption of Lemma 3.2, we have:
\[\int _{\theta }^{\theta +\lambda \varepsilon }\, \, E\left(\bar{L}\right)(\, t,\xi \, ;\, \dot{\bar{x}}(\cdot ))\, dt =\varepsilon \lambda \, E\left(\bar{L}\right)(\, \theta _{+} ,\xi \, ;\, \dot{\bar{x}}(\cdot ))+\frac{\varepsilon ^{2} \lambda ^{2} }{2} \frac{d}{dt} E\left(\bar{L}\right)(\, \theta _{+} ,\xi \, ;\, \dot{\bar{x}}(\cdot ))+o(\varepsilon ^{2} ),\] 
\[\int _{\theta +\lambda \varepsilon }^{\theta +\varepsilon }\, \, E\left(\bar{L}\right)\left(\, t,\frac{\lambda }{\lambda -1} \xi \, ;\, \dot{\bar{x}}(\cdot )\right)\, dt =\varepsilon (1-\lambda )\, E\left(\bar{L}\right)\left(\, \theta _{+} ,\frac{\lambda }{\lambda -1} \xi \, ;\, \dot{\bar{x}}(\cdot )\right)+\] 
\[+\frac{1}{2} \varepsilon ^{2} (1-\lambda ^{2} )\frac{d}{dt} E\left(\bar{L}\right)\left(\, \theta _{+} ,\frac{\lambda }{\lambda -1} \xi \, ;\, \dot{\bar{x}}(\cdot )\right)+o(\varepsilon ^{2} ).\] 

From here, considering \eqref{GrindEQ__1_5_}, we obtain:
\[\int _{\theta }^{\theta +\lambda \varepsilon }\, \, E\left(\bar{L}\right)(\, t,\xi \, ;\, \dot{\bar{x}}(\cdot ))\, dt +\int _{\theta +\lambda \varepsilon }^{\theta +\varepsilon }\, \, E\left(\bar{L}\right)\left(\, t,\frac{\lambda }{\lambda -1} \xi \, ;\, \dot{\bar{x}}(\cdot )\right)\, dt =\] 
\begin{equation} \label{GrindEQ__3_11_} 
=\varepsilon Q_{1} \left(\bar{L}\right)(\, \theta _{+} ,\lambda ,\xi \, ;\, \dot{\bar{x}}(\cdot ))+\frac{\varepsilon ^{2} }{2} \frac{d}{dt} Q_{2} \left(\bar{L}\right)(\, \theta _{+} ,\lambda ,\xi \, ;\, \dot{\bar{x}}(\cdot ))+o(\varepsilon ^{2} ).             
\end{equation} 

Quite similarly we have
\[\int _{\theta }^{\theta +\lambda \varepsilon }\, \, E\left(\bar{L}\right)(\, t+h,\xi \, ;\, \dot{\bar{y}}(\cdot ))\, dt +\int _{\theta +\lambda \varepsilon }^{\theta +\varepsilon }\, \, E\left(\bar{L}\right)\left(\, t+h,\frac{\lambda }{\lambda -1} \xi \, ;\, \dot{\bar{y}}(\cdot )\right)\, dt =\] 
\begin{equation} \label{GrindEQ__3_12_} 
=\varepsilon Q_{1} \left(\bar{L}\right)(\, (\theta +h)_{+} ,\lambda ,\xi \, ;\, \dot{\bar{y}}(\cdot ))+\frac{1}{2} \varepsilon ^{2} \frac{d}{dt} Q_{2} \left(\bar{L}\right)(\, (\theta +h)_{+} ,\lambda ,\xi \, ;\, \dot{\bar{y}}(\cdot ))+o(\varepsilon ^{2} ).   
\end{equation} 

Therefore, taking into accounts \eqref{GrindEQ__3_11_} and \eqref{GrindEQ__3_12_} in \eqref{GrindEQ__3_9_} we obtain the increment formula \eqref{GrindEQ__3_10_}. Lemma 3.2 is proven.

\section{Formulas for the increment of the functional in the problem (1.1),(1.2) on variations (2.3)}

Considering variation \eqref{GrindEQ__2_3_}, quite similarly to Lemma 3.1, the following Lemma is proven.

\textbf{Lemma 4.1.} Let functions $L(\cdot )$ and  $\varphi \left(\cdot \right)$ be continuously differentiable in the totality of their arguments and an admissible function $\bar{x}(\cdot )$ is an extremal of the problem \eqref{GrindEQ__1_1_}, \eqref{GrindEQ__1_2_}, i.e. solution of equation \eqref{GrindEQ__1_7_}. In addition, let the derivative function  $\dot{x}(\cdot )$ be continuous on the left at the points $\theta -h,\, \theta $ and $\theta +h$, where  $\theta \in (t_0,t_1-h]$. Then for the increment 
\begin{equation} \label{GrindEQ__4_1_} 
S(x^{(-)}(\cdot; \vartheta,\varepsilon);\overline{x}(\cdot ))-S(\overline{x}(\cdot ))=:\Delta^{(-)}S(q^{(-)}(\cdot; \vartheta,\varepsilon);\overline{x}(\cdot )) 
\end{equation} 
there is a representation of the form:
\[
\Delta^{(-)}S(q^{(-)}(\cdot ;\vartheta ,\varepsilon );\overline{x}(\cdot
))=\int_{\theta -\lambda \varepsilon }^{\theta }{\left[ E\left( \bar{L}%
\right) (t,\xi ;\dot{\bar{x}}(\cdot ))+E\left( \bar{L}\right) (t+h,\xi ;\dot{%
\bar{y}}(\cdot ))\right] dt}+
\]%
\[
+\int_{\theta -\varepsilon }^{\theta -\lambda \varepsilon }\,\left[ E\left( 
\bar{L}\right) \left( t,\frac{\lambda }{\lambda -1}\xi ;\,\dot{\bar{x}}%
(\cdot )\right) +E\left( \bar{L}\right) \left( t+h,\frac{\lambda }{\lambda -1
}\xi ;\,\dot{\bar{y}}(\cdot )\right) \right] dt\,-
\]%
\begin{equation}\label{GrindEQ__4_2_}
-\frac{1}{2}\varepsilon ^{2}\lambda \,\left[ M\left( \bar{L}_{x}^{{}}\right)
\left( \theta _{-},\lambda ,\xi ;\,\dot{\bar{x}}(\cdot )\right) +M\left( 
\bar{L}_{y}^{{}}\right) \left( (\theta +h)_{-},\lambda ,\xi ;\,\dot{\bar{y}}%
(\cdot )\right) \right] +o(\varepsilon ^{2}),  
\end{equation}
where functions $x^{(-)} \left(\cdot  ,\, \vartheta ,\varepsilon \right)$ and $q^{(-)} \left(\cdot  ,\, \vartheta \, \, ,\varepsilon \right)$ are determined by \eqref{GrindEQ__2_3_} and \eqref{GrindEQ__2_4_}, in addition $\vartheta=(\theta,\lambda,\xi)\in (t_0,t_1-h]\times (0,1)\times R^n$, $\varepsilon \in (0,\tilde{\varepsilon })$, $\tilde{\varepsilon }=\min \{ h,\theta -t_{0} \} $ and functions $E\left(\bar{L}\right)\left(\cdot \, \, ;\, \dot{\bar{x}}(\cdot )\right)$, $E\left(\bar{L}\right)\left(\cdot \, \, \, ;\, \dot{\bar{y}}(\cdot )\right)$, $M\left(\bar{L}_{x} \right)\left(\cdot \, ;\, \dot{\bar{x}}(\cdot )\right)$ and $M\left(\bar{L}_{y} \right)\left(\cdot \, \, ,\, \cdot \, \, ;\, \dot{\bar{y}}(\cdot )\right)$ are determined by \eqref{GrindEQ__1_4_} and \eqref{GrindEQ__1_6_} taking into account \eqref{GrindEQ__1_3_}.

\textbf{Proof.} Using the Taylor's formula with consideration of \eqref{GrindEQ__1_3_}, \eqref{GrindEQ__1_4_}, \eqref{GrindEQ__2_3_} and Lemma 2.2, we can easily state that for the increment \eqref{GrindEQ__4_1_}, the equality of the following form holds:
\[
\Delta ^{(-)}S(q^{(-)}(\cdot ;\vartheta ,\varepsilon );\overline{x}(\cdot
))=\Delta ^{(-)}S_{1}(q^{(-)}(\cdot ;\vartheta ,\varepsilon );\overline{x}%
(\cdot ))+
\]%
\begin{equation}\label{GrindEQ__4_3_}
+\Delta ^{(-)}S_{2}(q^{(-)}(\cdot ;\vartheta ,\varepsilon );\overline{x}%
(\cdot ))-\delta S(q^{(-)}(\cdot ;\vartheta ,\varepsilon );\overline{x}(\cdot )).
\end{equation}%
Here 
\[
\Delta ^{(-)}S_{1}(q^{(-)}(\cdot ;\vartheta ,\varepsilon );\overline{x}%
(\cdot ))=\int_{\theta -\lambda \varepsilon }^{\theta }{\left\{ \left[ \bar{L%
}(t,\xi ;\dot{\bar{x}}(\cdot ))-\bar{L}(t)\right] +\left[ \bar{L}(t+h,\xi ;%
\dot{\bar{y}}(\cdot ))-\bar{L}(t+h)\right] \right\} dt}+
\]%
\begin{equation}\label{GrindEQ__4_4_}
+\int_{\theta -\varepsilon }^{\theta -\lambda \varepsilon }\,\left\{ \,\left[
\bar{L}\left( t,\frac{\lambda }{\lambda -1}\xi \,;\dot{\bar{x}}(\cdot
)\right) -\bar{L}(\,t)\right] +\,\left[ \bar{L}\,\left( t+h,\frac{\lambda }{%
\lambda -1}\xi \,;\dot{\bar{y}}(\cdot )\right) -\bar{L}(\,t+h)\right]
\right\} \,dt,
\end{equation}%
\[
\Delta ^{(-)}S_{2}(q^{(-)}(\cdot ;\vartheta ,\varepsilon );\overline{x}%
(\cdot ))=\int_{\theta -\lambda \varepsilon }^{\theta }{\left\{ \left[ {\bar{%
L}}_{x}^{T}(t,\xi ;\dot{\bar{x}}(\cdot ))+{\bar{L}}_{y}^{T}(t+h,\xi ;\dot{%
\bar{y}}(\cdot ))\right] (t-\theta )\xi +o(t-\theta )\right\} dt}+
\]%
\[
+\frac{\lambda }{\lambda -1}\int_{\theta -\varepsilon }^{\theta -\lambda
\varepsilon }\,\left\{ \,\left[ \bar{L}_{x}^{{}}\left( t,\frac{\lambda }{%
\lambda -1}\xi \,;\dot{\bar{x}}(\cdot )\right) +\bar{L}_{y}^{{}}\,\left( t+h,%
\frac{\lambda }{\lambda -1}\xi \,;\dot{\bar{y}}(\cdot )\right) \right] \xi
(t-\theta +\varepsilon )+o(t-\theta +\varepsilon )\right\} \,dt=
\]%
\[
 =-\frac{1}{2}\lambda ^{2}\varepsilon ^{2}\left[ {%
\bar{L}}_{x}^{T}(\theta _{-},\xi ;\dot{\bar{x}}(\cdot ))+{\bar{L}}%
_{y}^{T}((\theta+h)_{-},\xi ;\dot{\bar{y}}(\cdot ))\right] \xi -
\]
\begin{equation}\label{GrindEQ__4_5_}
-\frac{1}{2}\varepsilon ^{2}\lambda (1-\lambda )\left[ {\bar{L}}%
_{x}^{T}\left( \theta _{-} ,\frac{\lambda }{\lambda -1}\xi;\dot{\bar{x}}
(\cdot )\right) +{\bar{L}}_{y}^{T}\left( {(\theta +h)}_{-},\frac{\lambda }{
\lambda -1}\xi ;\dot{\bar{y}}(\cdot )\right) \right] \xi +o(\varepsilon^{2}),
\end{equation}%
In addition, from \eqref{GrindEQ__2_8_} it follows that: 
\[
\delta S(q^{(-)}(\cdot ;\vartheta ,\varepsilon );\bar{x}(\cdot ))=\int_{\theta
-\lambda \varepsilon }^{\theta }{\left[ {\bar{L}}_{\dot{x}}^{T}(t)+{\bar{L}}%
_{\dot{y}}^{T}(t+h)\right] \xi dt}+
\]%
\[
 +\frac{\lambda }{\lambda -1}\int_{\theta -\varepsilon
}^{\theta -\lambda \varepsilon }{\left[ {\bar{L}}_{\dot{x}}^{T}\left(
t\right) +{\bar{L}}_{\dot{y}}^{T}\left( t+h\right) \right] \xi dt-\frac{1}{2}%
\lambda ^{2}\varepsilon ^{2}}\left[ {\bar{L}}_{x}^{T}\left( \theta
_{-}\right) +{\bar{L}}_{y}^{T}\left( {(\theta +h)}_{-}\right) \right] \xi -
\]
\begin{equation}\label{GrindEQ__4_6_}
-\frac{1}{2}\varepsilon ^{2}\lambda
(1-\lambda )\left[ {\bar{L}}_{x}^{T}\left( \theta _{-}\right) +{\bar{L}}%
_{y}^{T}\left( {(\theta +h)}_{-}\right) \right] \xi +o(\varepsilon ^{2}),
\end{equation}
where $q^{(-)} (\, \cdot \, \, ;\vartheta ,\varepsilon )$ is determined by \eqref{GrindEQ__2_4_}.

Similar to \eqref{GrindEQ__3_9_} from \eqref{GrindEQ__4_3_} considering \eqref{GrindEQ__1_3_}, \eqref{GrindEQ__1_4_}, \eqref{GrindEQ__1_6_} and \eqref{GrindEQ__4_4_} - \eqref{GrindEQ__4_6_} it is not difficult to obtain the validity of the expansion formula \eqref{GrindEQ__4_2_}. Lemma 4.1 is proven.

Similar to Lemma 3.2, using Lemma 4.1 and assuming additional smoothness conditions, the following statement is proven.

\textbf{Lemma 4.2}. Let functions $L(\cdot )$, $L_{\dot{x}} (\cdot )$ and $L_{\dot{y}} (\cdot )$ be continuously differentiable in the totality of their arguments and the function $\varphi\left(\cdot \right)-$ twice continuously differentiable. Moreover, let the admissible function $\bar{x}(\cdot )$ be an extremal of the problem \eqref{GrindEQ__1_1_}, \eqref{GrindEQ__1_2_} and it is twice continuously differentiable on the left in the semi-neighborhoods of the points $\theta -h,\, \theta $ $\theta +h$, where $\theta\in (t_0,t_1-h]{\rm .\ }$ Then for sufficiently small values $\varepsilon >0$ there is a representation of the form:
\[
\Delta ^{(-)}S(q^{(-)}(\cdot ;\vartheta ,\varepsilon );\overline{x}(\cdot
))=\varepsilon \left[ Q_{1}\left( \bar{L}\right) (\theta _{-},\lambda ,\xi ;%
\dot{\bar{x}}(\cdot ))+Q_{1}\left( \bar{L}\right)((\theta +h)_{-},\lambda ,\xi ;\dot{\bar{y}}%
(\cdot ))\right] -
\]%
\[
+\frac{1}{2}\varepsilon ^{2}\left\{ \lambda \,\left[ M\left( \bar{L}%
_{x}^{{}}\right) \left( \theta _{-},\lambda ,\xi ;\,\dot{\bar{x}}(\cdot
)\right) +M\left( \bar{L}_{y}^{{}}\right) \left( (\theta +h)_{-},\lambda
,\xi ;\,\dot{\bar{y}}(\cdot )\right) \right] \right. +
\]%
\begin{equation}
\left. +\frac{d}{dt}\left[ Q_{2}\left( \bar{L}\right) (\theta _{-},\lambda
,\xi ;\dot{\bar{x}}(\cdot ))+Q_{2}\left( \bar{L}\right)((\theta +h)_{-},\lambda ,\xi ;\dot{\bar{y}%
}(\cdot ))\right] \right\} +o(\varepsilon ^{2}),  \label{GrindEQ__4_7_}
\end{equation}
where $Q_{k} \left(\bar{L}\right)\left(\, \cdot \, ;\, \dot{\bar{x}}(\cdot )\right)$ and $Q_{k} \left(\bar{L}\right)\left(\, \cdot \, ;\, \dot{\bar{y}}(\cdot )\right),\, k=1,2$ ,$M\left(\bar{L}_{x}^{} \right)\left(\, \cdot \, ;\, \dot{\bar{x}}(\cdot )\right)$, $M\left(\bar{L}_{y}^{} \right)\left(\, \cdot \, ;\, \dot{\bar{y}}(\cdot )\right)$ are determined by \eqref{GrindEQ__1_5_} and \eqref{GrindEQ__1_6_}.

\textbf{Proof.} Due to the assumption of Lemma 4.2, it is asserted that the expansion \eqref{GrindEQ__4_2_} holds. Therefore, it is clear that to prove Lemma 4.2 it is enough to calculate the integrals in \eqref{GrindEQ__4_2_}. The calculation of these integrals is carried out similarly to \eqref{GrindEQ__3_11_} and \eqref{GrindEQ__3_12_}. Specifically, considering \eqref{GrindEQ__1_5_} and the assumption of Lemma 4.2, we have:
\[\int _{\theta -\lambda \varepsilon }^{\theta }\, \, E\left(\bar{L}\right)(\, t,\xi \, ;\, \dot{\bar{x}}(\cdot ))\, dt +\int _{\theta -\varepsilon }^{\theta -\lambda \varepsilon }\, E\left(\bar{L}\right)\left(t,\frac{\lambda }{\lambda -1} \xi ;\, \dot{\bar{x}}(\cdot )\right)dt\,  =\] 
\[=\varepsilon Q_{1} \left(\bar{L}_{}^{} \right)\left(\theta _{-} ,\lambda ,\xi ;\, \dot{\bar{x}}(\cdot )\right)-\frac{\varepsilon ^{2} }{2} \frac{d}{dt} Q_{2} \left(\bar{L}_{}^{} \right)\left(\theta _{-} ,\lambda ,\xi ;\, \dot{\bar{x}}(\cdot )\right)+o(\varepsilon ^{2} ),\] 
\[\int _{\theta -\lambda \varepsilon }^{\theta }\, \, E\left(\bar{L}\right)(\, t+h,\xi \, ;\, \dot{\bar{y}}(\cdot ))\, dt +\int _{\theta -\varepsilon }^{\theta -\lambda \varepsilon }\, E\left(\bar{L}\right)\left(t+h,\frac{\lambda }{\lambda -1} \xi ;\, \dot{\bar{y}}(\cdot )\right)dt\,  =\] 
\[=\varepsilon Q_{1} \left(\bar{L}_{}^{} \right)\left((\theta +h)_{-} ,\lambda ,\xi ;\, \dot{\bar{y}}(\cdot )\right)-\frac{1}{2} \varepsilon ^{2} \frac{d}{dt} Q_{2} \left(\bar{L}_{}^{} \right)\left((\theta +h)_{-} ,\lambda ,\xi ;\, \dot{\bar{y}}(\cdot )\right)+o(\varepsilon ^{2} ).\]

First, we summarize the last equalities and then take them into account in \eqref{GrindEQ__4_2_}. As a result, we obtain the increment formula 4.7. Therefore, Lemma 4.2 is proven.

\section{ Necessary conditions for a minimum in the presence of degeneracy on an interval}

At this point, using Lemma 3.1 and 4.1, i.e. based on the increment formulas \eqref{GrindEQ__3_9_} and \eqref{GrindEQ__4_2_}, the following theorem is proved.

\textbf{Theorem 5.1.} Let functions $L(\cdot )$ and  $\varphi \left(\cdot \right)$ be continuously differentiable with respect to their arguments and an admissible function $\bar{x}(\cdot )$ is an extremal of the problem \eqref{GrindEQ__1_1_}, \eqref{GrindEQ__1_2_}, then along it for vectors $\eta \in R^{n} \backslash \{ 0\} $ and $\frac{\bar{\lambda }}{\bar{\lambda }-1} \eta ,\, \, \, \bar{\lambda }\in (0,\, 1)$ Weierstrass condition at any point of the interval $\left(\bar{t}_{0} ,\bar{t}_{1} \right)\subset [t_{0} ,t_{1} -h]$ degenerates, i.e. the equalities hold:
\begin{eqnarray} \label{GrindEQ__5_1_} 
E\left( \bar{L}\right) (t,\eta ;\dot{\bar{x}}(\cdot ))+E\left( \bar{L}%
\right) (t+h,\eta ;\dot{\bar{y}}(\cdot )) &=&0,\forall t\in ({\bar{t}}_{0},{%
\bar{t}}_{1}), \notag \\
E\left( \bar{L}\right) \left( t,\frac{\overline{\lambda }}{\overline{\lambda 
}-1}\eta ;\dot{\bar{x}}(\cdot )\right) +E\left( \bar{L}\right) \left( t,%
\frac{\overline{\lambda }}{\overline{\lambda }-1}\eta ;\dot{\bar{y}}(\cdot
)\right)  &=&0,\forall t\in ({\bar{t}}_{0},{\bar{t}}_{1}).
\end{eqnarray}

In addition, let the function $\bar{x}(\cdot )$ be continuously differentiable in the intervals \\ $\left(\bar{t}_{0} -h,\bar{t}_{1} -h\right), (\bar{t}_{0} ,\bar{t}_{1} )$ and $\left(\bar{t}_{0} +h,\bar{t}_{1} +h\right)$. Then:

(i)    if the function $\bar{x}(\cdot )$ is a strong local minimum in the problem \eqref{GrindEQ__1_1_}, \eqref{GrindEQ__1_2_}, then the equality holds:
\begin{equation} \label{GrindEQ__5_2_} 
M\left(\bar{L}_{x} \right)(t,\bar{\lambda },\eta ;\, \, \dot{\bar{x}}(\cdot ))+M\left(\bar{L}_{y} \right)(t+h,\bar{\lambda },\eta ;\, \, \dot{\bar{y}}(\cdot ))=0,\, \, \, \, \, \, \forall \, \, t\in (\bar{t}_{0} ,\bar{t}_{1} )\, ,           
\end{equation} 

where functions $M\left(\bar{L}_{x} \right)\left(\cdot \, \, \, ;\, \dot{\bar{x}}(\cdot )\right)$ and $M\left(\bar{L}_{y} \right)\left(\cdot \, \, ;\, \dot{\bar{y}}(\cdot )\right)$ are determined by \eqref{GrindEQ__1_6_};

(ii)   if the function  $\bar{x}(\cdot )$ is a weak local minimum in the problem \eqref{GrindEQ__1_1_}, \eqref{GrindEQ__1_2_}, then there exist  a number $\delta >0$, at which for each $ \left(\overline{\lambda },\eta, \frac{\overline{\lambda }}{\overline{\lambda }-1}\eta \right) \in  (0,1)\times B_{\delta } (0) \times B_{\delta } (0)$, satisfying condition \eqref{GrindEQ__5_1_}, equality \eqref{GrindEQ__5_2_} holds, where $B_{\delta } (0)$ -  is a closed ball of radius $\delta $ centered at  the point $0\in R^n$.

\textbf{Proof.}  Firstly, let's prove part (i) of Theorem 5.1. Due to the assumption of Theorem 5.1, we assert: firstly, Lemmas 3.1 and 4.1 are true, i.e. expansion formulas \eqref{GrindEQ__3_9_} and \eqref{GrindEQ__4_2_} are valid; secondly, for each $\vartheta =(\theta ,\lambda ,\xi )\in (\bar{t}_{0} ,\bar{t}_{1} )\times (0,\, 1)\times R^{n} $ and for sufficiently small $\varepsilon >0$ ones the inequalities hold:
\begin{equation} \label{GrindEQ__5_3_} 
\Delta ^{(+)} S_{} (q^{(+)} (\, \cdot \, ;\, \vartheta ,\varepsilon );{\bar{x}}(\cdot ))\ge 0,\, \, \, \, \, \, \, \Delta ^{(-)} S_{} (q^{(-)} (\, \cdot \, ;\, \vartheta ,\varepsilon );{\bar{x}}(\cdot ))\, \, \ge 0.           
\end{equation}

Let's put $\vartheta =\overline{\vartheta}:=(\theta ,\overline{\lambda },\eta
)\in (\overline{t}_{0},\overline{t}_{1})\times (0,1)\times R^{n}\backslash
\{0\}$,\textit{\ } where $\theta $ is arbitrary fixed point. Obviously, you can
choose a number $\varepsilon^{\ast }>0$ so that there is a shutdown $[\theta
-\varepsilon ^{\ast },\theta +\varepsilon ^{\ast }]\subset (\overline{t}_{0},%
\overline{t}_{1})$. Let's consider this in the formulas \eqref{GrindEQ__3_9_}
and \eqref{GrindEQ__4_2_}. Then, by virtue of \eqref{GrindEQ__5_1_}, for
sufficiently small $\varepsilon \in (0,\varepsilon ^{\ast }]$ , the first
two terms in \eqref{GrindEQ__3_9_}, as well as in \eqref{GrindEQ__4_2_},
i.e. all integral terms in \eqref{GrindEQ__3_9_} and \eqref{GrindEQ__4_2_}
become zero. As a result, the expansion formulas \eqref{GrindEQ__3_9_} and %
\eqref{GrindEQ__4_2_} considering \eqref{GrindEQ__5_3_} take the following
form:

\[
\Delta ^{(+)}S(q^{(+)}(\cdot ;\overline{\vartheta },\varepsilon )
;\overline{x}(\cdot ))=
\]
\begin{equation}\label{GrindEQ__5_4_}
=\frac{1}{2}\varepsilon ^{2}\overline{\lambda }\left[ M\left( \overline{{L}}%
_{x}\right) (\theta _{+},\overline{\lambda },\eta ;\,\dot{\bar{x}}(\cdot)+M({\bar{L}}%
_{y})((\theta +h)_{+},\overline{\lambda },\eta ;\dot{\bar{y}}(\cdot ))\right]
+o(\varepsilon ^{2})\ge 0,  
\end{equation}
\[
\Delta ^{(-)}S(q^{(-)}(\cdot ;\overline{\vartheta },\varepsilon );\bar{x}%
(\cdot ))=
\]
\begin{equation}\label{GrindEQ__5_5_}
=\frac{1}{2}\varepsilon ^{2}\overline{\lambda }\left[ M\left( {\bar{L}}%
_{x}\right) (\theta _{-},,\overline{\lambda },\eta ;\dot{\bar{x}}(\cdot ))+M(%
{\bar{L}}_{y})((\theta +h)_{-},\overline{\lambda },\eta ;\dot{\bar{y}}(\cdot
))\right] -o(\varepsilon ^{2})\le 0,  
\end{equation}
where $\varepsilon \in (0,\varepsilon ^{*} ]$ is quite small number and functions $M\left(\bar{L}_{x} \right)\left(\cdot \, \right)$ and $M\left(\bar{L}_{y} \right)\left(\cdot \, \right)$ are determined by \eqref{GrindEQ__1_6_}.

By virtue of the assumption of smoothness of the function $L(\cdot )$, $\varphi\left(\cdot \right)$ and $\bar{x}(\cdot )$, and considering \eqref{GrindEQ__1_6_}, we can easily confirm that the coefficients $\varepsilon ^{2} $ in inequalities \eqref{GrindEQ__5_4_} and \eqref{GrindEQ__5_5_} coincide. Therefore, we obtain:
\[
\mathop{\lim }\limits_{\varepsilon \rightarrow +0}\,\,\frac{1}{\varepsilon
^{2}}\,\,\Delta ^{(+)}S(q^{(+)}(\,\cdot \,;\,\bar{\vartheta},\varepsilon );%
\bar{x}(\cdot ))\,\,=\mathop{\lim }\limits_{\varepsilon \rightarrow +0}\,\,%
\frac{1}{\varepsilon ^{2}}\Delta ^{(-)}S(q^{(-)}(\,\cdot \,;\,\bar{\vartheta}%
,\varepsilon );\bar{x}(\cdot ))=
\]%
\[
=\frac{1}{2}\overline{\lambda }D(\theta ,\overline{\lambda },\eta ;\bar{x}%
(\cdot )),
\]
where $D(\theta ,\bar{\lambda },\eta ;\bar{x}(\cdot ))\, \, =M\left(\bar{L}_{x} \right)(\, \theta ,\bar{\lambda },\eta \, ;\, \dot{\bar{x}}(\cdot ))+M(\bar{L}_{y} )(\, \theta ,\bar{\lambda },\eta \, ;\, \dot{\bar{y}}(\cdot ))$.

It is clear that due to the inequalities \eqref{GrindEQ__5_4_} and \eqref{GrindEQ__5_5_} we have $D(\theta ,\bar{\lambda },\eta ;\bar{x}(\cdot ))\ge 0$ and $D(\theta ,\bar{\lambda },\eta ;\bar{x}(\cdot ))\le 0$. Hence, considering the arbitrariness of  $\theta \in (\bar{t}_{0} ,\bar{t}_{1} )$, the proof of equality \eqref{GrindEQ__5_2_} follows. Consequently, part (i) of Theorem 5.1 is proven.


Now we prove part (ii) of Theorem 5.1. Let us consider the increment formulas \eqref{GrindEQ__5_4_} and \eqref{GrindEQ__5_5_} obtained for $\Delta ^{(+)} S(q^{(+)} (\, \cdot \, ;\, \bar{\vartheta },\varepsilon );\bar{x}(\cdot ))$ and $\Delta ^{(-)} S(q^{(-)} (\, \cdot \, ;\, \bar{\vartheta },\varepsilon );\bar{x}(\cdot ))$. Here

\[\Delta^{(+)}S(q^{(+)}(\cdot ;\bar{\vartheta},\varepsilon);\bar{x}(\cdot ))=S(x^{(+)}(\cdot ;\bar{\vartheta},\varepsilon))-S(\bar{x}(\cdot ))\]
and  
\[\Delta^{(-)}S(q^{(-)}(\cdot ;\bar{\vartheta},\varepsilon);\bar{x}(\cdot ))=S(x^{(-)}(\cdot ;\bar{\vartheta},\varepsilon))-S(\bar{x}(\cdot )),\] 
where the functions $x^{(+)} (\, \cdot \, ;\, \bar{\vartheta },\varepsilon )$ and $x^{(-)} (\, \cdot \, ;\, \bar{\vartheta },\varepsilon )$  for  $\vartheta =\bar{\vartheta } $ are determined by \eqref{GrindEQ__2_1_}, \eqref{GrindEQ__2_2_} and \eqref{GrindEQ__2_3_}, \eqref{GrindEQ__2_4_}, respectively.

Further, by virtue of the definition of the function $x^{(+)} (\cdot \, ;\, \bar{\vartheta },\varepsilon )$ and $x^{(-)} (\, \cdot \, ;\, \bar{\vartheta },\varepsilon )$ considering \eqref{GrindEQ__2_5_} and \eqref{GrindEQ__2_6_} for all $\varepsilon\in (0,\overline{\varepsilon})\cap \ \ (0,1]$ the following estimates are valid:
\begin{equation}\label{GrindEQ__5_6_}
{\left\Vert x^{(+)}(\cdot ;\bar{\vartheta},\varepsilon )-\bar{x}(\cdot
)\right\Vert }_{C\left( \hat{I},R^{n}\right) }={\left\Vert q^{(+)}(\cdot ;%
\bar{\vartheta},\varepsilon )\right\Vert }_{C\left( \hat{I}%
,R^{n}\right) }\le \left\Vert \eta \right\Vert
\end{equation}%
\begin{equation}\label{GrindEQ__5_7_}
{\left\Vert {\dot{x}}^{(+)}(\cdot ;\bar{\vartheta},\varepsilon )-\dot{%
\overline{x}}(\cdot )\right\Vert }_{L_{\infty}\left( \hat{I},R^{n}\right) }={%
\left\Vert {\dot{q}}^{(+)}(\cdot ;\bar{\vartheta},\varepsilon )\right\Vert }%
_{L_{\infty}\left( \hat{I},R^{n}\right) }=\mathrm{max}\left\{{\left\Vert \eta
\right\Vert }_{R^{n}},\frac{\overline{\lambda }}{1-\overline{\lambda }}{%
\left\Vert \eta \right\Vert }_{R^{n}}\right \}.
\end{equation}

Similar estimates are valid for $\left\| x^{(-)} (\cdot; \bar{\vartheta },\varepsilon )-\bar{x}(\cdot )\right\| _{C\left(\hat{I},R^{n} \right)} $ and  ${\left\|{\dot{x}}^{(-)}(\cdot ;\bar{\vartheta},\varepsilon)-\dot{\bar{x}}(\cdot )\right\|}_{L_\infty \left(\hat{I},R^n\right)}$ at        $\varepsilon\in (0,\widetilde{\varepsilon })\cap \ \ (0,1]$.

Since, according to conditions part (ii)  of Theorem 5.1, the admissible function is a weak local minimum in the problem \eqref{GrindEQ__1_1_}, \eqref{GrindEQ__1_2_} (for certainty with a neighborhood $\hat{\delta }$), then considering estimates \eqref{GrindEQ__5_6_} and \eqref{GrindEQ__5_7_} for each point $\left(\eta ,(\bar{\lambda }-1)^{-1} \bar{\lambda }\eta ,\, \bar{\lambda }\right)\in B_{\hat{\delta }} (0)\times B_{\hat{\delta }} (0)\times (0,1)$ and   
$\varepsilon\in \left(0,\hat{\varepsilon}\right),\hat{\varepsilon}=min\{\overline{\varepsilon},\tilde{\varepsilon},1\}$, satisfying condition \eqref{GrindEQ__5_1_}, inequalities \eqref{GrindEQ__5_4_} and \eqref{GrindEQ__5_5_} hold. Based on the last derivation, by choosing $\delta =\hat{\delta }$ and considering the arbitrariness $\theta \in (\bar{t}_{0} ,\bar{t}_{1} )$, it is not difficult to obtain a proof of part (ii) of Theorem 5.1. Therefore, Theorem 5.1 is proven.

\section{ Necessary minimum conditions in the presence of degeneracy at the point.}

In this section, unlike Section 5, we investigate the minimum of the extremal in problem \eqref{GrindEQ__1_1_}, \eqref{GrindEQ__1_2_} in the presence of degeneracy at one point. Specifically, with additional assumptions of smoothness, using Lemmas 3.2 and 4.2, i.e., based on formulas \eqref{GrindEQ__3_10_} and \eqref{GrindEQ__4_7_}, the following statement is proven.

\textbf{Theorem 6.1.} Let the functions $L(\cdot )$, $L_{\dot{x}} (\cdot )$ and $L_{\dot{y}} (\cdot )$ be continuously differentiable with respect to their arguments and the function  $\varphi\left(\cdot \right)-$ twice continuously differentiable. In addition, let the admissible function $\bar{x}(\cdot )$ be a strong local minimum of the problem \eqref{GrindEQ__1_1_}, (1, 2). Then:

(i) if $\theta \in [t_{0} ,t_{1} -h)$  ($\theta \in (t_{0} ,t_{1} -h]$) and the function $\bar{x}(\cdot )$ is twice continuously differentiable in the right (left) semi-vicinity of the points $\theta -h,\, \theta $ and $\theta +h$, in addition, along it, for the number  $\bar{\lambda }\in (0,\, 1)$, as well as for the vectors $\eta \ne 0$ and $(\bar{\lambda}-1)^{-1}\bar{\lambda}\eta$  Weierstrass condition \eqref{GrindEQ__1_9_} degenerates on the right (left) at the point $\theta $,  i.e. the equalities hold:
\[
E\left( \bar{L}\right) (\theta _{+},\eta ;\,\,\dot{\bar{x}}(\cdot ))+E\left( 
\bar{L}\right) ((\theta +h)_{+},\eta ;\,\,\dot{\bar{y}}(\cdot ))=
\]%
\begin{equation}
=E\left( \bar{L}\right) \left( \theta _{+},(\overline{\lambda }-1)^{-1}%
\overline{\lambda }\eta ;\dot{\bar{x}}(\cdot )\right) +E\left( \bar{L}%
\right) \left( (\theta +h)_{+},(\overline{\lambda }-1)^{-1}\overline{\lambda 
}\eta ;\dot{\bar{y}}(\cdot )\right) =0  \label{GrindEQ__6_1_}
\end{equation}%
\[
\left( E\left( \bar{L}\right) (\theta _{-},\eta ;\,\,\dot{\bar{x}}(\cdot
))+E\left( \bar{L}\right) ((\theta +h)_{-},\eta ;\,\,\dot{\bar{y}}(\cdot
))=\right. 
\]%
\begin{equation}
\left. =E\left( \bar{L}\right) \left( \theta _{-},(\overline{\lambda }%
-1)^{-1}\overline{\lambda }\eta ;\dot{\bar{x}}(\cdot )\right) +E\left( \bar{L%
}\right) \left( (\theta +h)_{-},(\overline{\lambda }-1)^{-1}\overline{%
\lambda }\eta ;\dot{\bar{y}}(\cdot )\right) =0\right) ,
\label{GrindEQ__6_2_}
\end{equation}%
then the following inequalities hold: 
\[
\bar{\lambda}\left[ M\left( \bar{L}_{x}\right) (\theta _{+},\bar{\lambda}%
,\eta ;\,\,\dot{\bar{x}}(\cdot ))+M\left( \bar{L}_{y}\right) ((\theta
+h)_{+},\bar{\lambda},\eta ;\,\,\dot{\bar{y}}(\cdot ))\right] +
\]%
\begin{equation}\label{GrindEQ__6_3_}
+\frac{d}{dt}\left[ Q_{2}\left( \bar{L}\right) \left( \theta _{+},\bar{%
\lambda},\eta ;\,\,\dot{\bar{x}}(\cdot )\right) +Q_{2}\left( \bar{L}\right)
\left( (\theta +h)_{+},\bar{\lambda},\eta ;\,\,\,\dot{\bar{y}}(\cdot
)\right) \right] \geq 0  
\end{equation}%
\[
\left( ^{^{{}}}\bar{\lambda}\left[ M\left( \bar{L}_{x}\right) (\theta _{-},%
\bar{\lambda},\eta ;\,\,\dot{\bar{x}}(\cdot ))+M\left( \bar{L}_{y}\right)
((\theta +h)_{-},\bar{\lambda},\eta ;\,\,\dot{\bar{y}}(\cdot ))\right]
^{^{{}}}\right. +
\]%
\begin{equation}\label{GrindEQ__6_4_}
\left. +\frac{d}{dt}\left[ Q_{2}\left( \bar{L}\right) \left( \theta _{-},%
\bar{\lambda},\eta ;\,\,\dot{\bar{x}}(\cdot )\right) +Q_{2}\left( \bar{L}%
\right) \left( (\theta +h)_{-},\bar{\lambda},\eta ;\,\,\,\dot{\bar{y}}(\cdot
)\right) \right] \leq 0^{^{{}}}\right) ,  
\end{equation}%
when $E\left( \bar{L}\right) (\cdot ),\,\,Q_{2}\left( \bar{L}\right) (\cdot
),\,\,M\left( \bar{L}_{x}\right) (\cdot )$ and $M\left( \bar{L}_{y}\right)
(\cdot )$are determined by \eqref{GrindEQ__1_4_} - \eqref{GrindEQ__1_6_};

(ii) if $\theta \in (t_{0} ,t_{1} -h)$ and the function $\bar{x}(\cdot )$ is twice continuously differentiable  in certain neighborhood of each points $\theta -h,\, \theta $ and $\theta +h$,  in addition, along it for the number $\bar{\lambda }\in (0,\, 1)$, as well as for the vectors  $\eta \ne 0$ and $(\bar{\lambda }-1)^{-1} \bar{\lambda }\eta $  the Weierstrass condition degenerates at the point $\theta $, i.e. the equalities hold
\[
E\left( \bar{L}\right) (\theta ,\eta ;\,\,\dot{\bar{x}}(\cdot ))+E\left( 
\bar{L}\right) (\theta +h,\eta ;\,\,\dot{\bar{y}}(\cdot )) =
\]
\begin{equation}\label{GrindEQ__6_5_}
E\left( \bar{L}\right) \left( \theta ,(\bar{\lambda}-1)^{-1}\bar{\lambda}%
\eta ;\dot{\bar{x}}(\cdot )\right) +E\left( \bar{L}\right) \left( \theta +h,(%
\bar{\lambda}-1)^{-1}\bar{\lambda}\eta ;\dot{\bar{y}}(\cdot )\right)  =0,
\end{equation} 
then the equality of the form holds
\begin{equation} \label{GrindEQ__6_6_} 
M\left(\bar{L}_{x} \right)(\theta ,\bar{\lambda },\eta ;\, \, \dot{\bar{x}}(\cdot ))+M\left(\bar{L}_{y} \right)(\theta +h,\bar{\lambda },\eta ;\, \, \dot{\bar{y}}(\cdot ))=0 .  
\end{equation} 

\textbf{Proof.} First, let us prove part (i) of Theorem 6.1. By virtue of the assumption of Theorem 6.1, Lemma 3.2 (Lemma 4.2) holds, i.e. the expansion \eqref{GrindEQ__3_10_} \eqref{GrindEQ__4_7_} is valid.

Put $\lambda =\bar{\lambda }$ and $\xi =\eta $ at \eqref{GrindEQ__3_10_} ((4.7)). Then by virtue of \eqref{GrindEQ__1_5_} and \eqref{GrindEQ__6_1_} ((6.2)) it is easy to obtain:
\begin{equation} \label{GrindEQ__6_7_} 
Q_{1} \left(\bar{L}\right)\left(\theta _{+} ,\bar{\lambda },\eta ;\, \, \dot{\bar{x}}(\cdot )\right)+Q_{1} \left(\bar{L}\right)\left((\theta +h)_{+} ,\bar{\lambda },\eta ;\, \, \, \dot{\bar{y}}(\cdot )\right)=0 
\end{equation} 
\begin{equation} \label{GrindEQ__6_8_} 
\left(\, Q_{1} \left(\bar{L}\right)\left(\theta _{-} ,\bar{\lambda },\eta ;\, \, \dot{\bar{x}}(\cdot )\right)+Q_{1} \left(\bar{L}\right)\left((\theta +h)_{-} ,\bar{\lambda },\eta ;\, \, \, \dot{\bar{y}}(\cdot )\right)=0\, \right).              
\end{equation} 
Since the function $\bar{x}(\cdot )$  is the strong local minimum of the problem \eqref{GrindEQ__1_1_}, \eqref{GrindEQ__1_2_}, from \eqref{GrindEQ__3_10_} (4.7) considering $\vartheta=\bar{\vartheta}=(\theta,\bar{\lambda},\eta)$ and the equation \eqref{GrindEQ__6_7_} ((6.8)) for sufficiently small $\varepsilon >0$  we have: 
\[
\varepsilon ^{-2}\Delta ^{(+)}S(q^{(+)}(\cdot ;\bar{\vartheta},\varepsilon );%
\bar{x}(\cdot ))-o(\varepsilon ^{2})/\varepsilon ^{2}\ge 0
\]%
\[
\left( \,\,\varepsilon ^{-2}\Delta ^{(-)}S(q^{(-)}(\,\cdot \,;\,\bar{%
\vartheta},\varepsilon );\bar{x}(\cdot ))-o(\varepsilon ^{2})/\varepsilon
^{2}\geq 0\,\,\right) .
\]

From here, passing to the limit  $\varepsilon \to +0$ considering \eqref{GrindEQ__3_10_} ((4.7)), we obtain the validity of inequality \eqref{GrindEQ__6_3_} ((6.4)). So, part (i) of Theorem 6.1 is proven. 

Now we prove part (ii) of Theorem 6.1. By virtue of assumption part (ii) of Theorem 6.1, it is easy to confirm that the left sides of inequalities \eqref{GrindEQ__6_3_} and \eqref{GrindEQ__6_4_} coincide. Therefore, the equality is valid:
\[\bar{\lambda }\left[M\left(\bar{L}_{x} \right)(\theta ,\bar{\lambda },\eta ;\, \, \dot{\bar{x}}(\cdot ))+M\left(\bar{L}_{y} \right)(\theta +h,\bar{\lambda },\eta ;\, \, \dot{\bar{y}}(\cdot ))\right]^{} +\] 
\begin{equation} \label{GrindEQ__6_9_} 
+\frac{d}{dt} \left[Q_{2} \left(\bar{L}\right)\left(\theta ,\bar{\lambda },\eta ;\, \, \dot{\bar{x}}(\cdot )\right)+Q_{2} \left(\bar{L}\right)\left(\theta +h,\bar{\lambda },\eta ;\, \, \, \dot{\bar{y}}(\cdot )\right)\right]=0.                 
\end{equation} 

Since the functions $Q_{2} \left(\bar{L}\right)\left(\, \cdot \, ,\bar{\lambda },\bar{\eta };\, \, \dot{\bar{x}}(\cdot )\right)$ and $Q_{2} \left(\bar{L}\right)\left(\, \cdot \, ,\bar{\lambda },\bar{\eta };\, \, \dot{\bar{y}}(\cdot )\right)$ are defined by \eqref{GrindEQ__1_5_} at $k=2$ , then, considering the assumptions about the smoothness of the function $L(\cdot )$, $L_{\dot{x}} (\cdot )$, $L_{\dot{y}} (\cdot )$, $\varphi(\cdot )$ and $\bar{x}(\cdot )$  we have:
\[
{\left. \frac{d}{dt}{\left[ Q_{2}\left( \bar{L}\right) \left( t,\bar{\lambda}%
,\eta ;\dot{\bar{x}}(\cdot )\right) +Q_{2}\left( \bar{L}\right) \left( t+h,%
\bar{\lambda},\eta ;\dot{\bar{y}}(\cdot )\right) \right] }\right\vert }%
_{t=\theta }=
\]%
\[
=\bar{\lambda}^{2}{\left. \frac{d}{dt}{\left[ E\left( \bar{L}\right) \left(
t,\eta ;\dot{\bar{x}}(\cdot )\right) +E\left( \bar{L}\right) \left( t+h,\eta
;\dot{\bar{y}}(\cdot )\right) \right] }\right\vert }_{t=\theta }+
\]%
\begin{equation}\label{GrindEQ__6_10_}
+\left( 1-\bar{\lambda}^{2}\right) {\left. \frac{d}{dt}{\left[ E\left( \bar{L%
}\right) \left( t,(\bar{\lambda}-1)^{-1}\bar{\lambda}\eta ;\dot{\bar{x}}%
(\cdot )\right) +E\left( \bar{L}\right) \left( t+h,(\bar{\lambda}-1)^{-1}%
\bar{\lambda}\eta ;\dot{\bar{y}}(\cdot )\right) \right] }\right\vert }%
_{t=\theta }
\end{equation}
Further, since the function $\bar{x}(\cdot )$  is a strong local minimum of the problem \eqref{GrindEQ__1_1_}, \eqref{GrindEQ__1_2_}, then, by virtue of Weierstrass condition \eqref{GrindEQ__1_9_} and assumption \eqref{GrindEQ__6_5_} of the functions

\[E\left(\bar{L}\right)\left(\cdot \, ,\eta ;\, \, \dot{\bar{x}}(\cdot )\right)+E\left(\bar{L}\right)\left(t+h,\eta ;\, \, \, \dot{\bar{y}}(\cdot )\right),\; t\in [t_{0} ,t_{1} -h]\]  
and  
\[E\left(\bar{L}\right)\left(t,(\bar{\lambda }-1)^{-1} \bar{\lambda }\eta \right.\left. \dot{\bar{x}}(\cdot )\right)+E\left(\bar{L}\right)\left(t+h,(\bar{\lambda }-1)^{-1} \bar{\lambda }\eta ;\, \, \, \dot{\bar{y}}(\cdot )\right),\;  t\in [t_{0} ,t_{1} -h]\]
for the variable $t$  reach the minimum at the point $\theta \in (t_{0} ,t_{1} -h)$. Then, taking into account the assumptions about the smoothness of the functions $L(\cdot )$, $L_{\dot{x}} (\cdot )$, $L_{\dot{y}} (\cdot )$, $\varphi(\cdot )$ and  $\bar{x}(\cdot )$ according to Fermat's theorem [9, page 15], we confirm that the derivatives of these functions at the point $\theta $ are equal zero, and therefore the left side \eqref{GrindEQ__6_10_} is equal  zero. 

 Therefore, from \eqref{GrindEQ__6_9_}, considering $\bar{\lambda }\in (0,\, 1)$  the last statement and inclusion, we obtain the proof of the equality \eqref{GrindEQ__6_6_}. Hence, part (ii) of theorem 6.1 is proven. Thus, theorem 6.1 is completely proven.

\textbf{Remark 6.1.} If the function $\bar{x}(\cdot )$ is a strong local minimum of the problem \eqref{GrindEQ__1_1_}, \eqref{GrindEQ__1_2_}, then along it from the equation:
\begin{equation} \label{GrindEQ__6_11_} 
Q_{1} \left(\bar{L}\right)\left(\hat{\theta },\bar{\lambda },\eta ;\, \, \dot{\bar{x}}(\cdot )\right)+Q_{1} \left(\bar{L}\right)\left(\hat{\theta }+h,\bar{\lambda },\eta ;\, \, \, \dot{\bar{y}}(\cdot )\right)=0 
\end{equation} 
considering \eqref{GrindEQ__1_5_} and \eqref{GrindEQ__1_9_} follows \eqref{GrindEQ__6_1_} and \eqref{GrindEQ__6_2_}, and vice versa, i.e. from \eqref{GrindEQ__6_1_} and \eqref{GrindEQ__6_2_} follows \eqref{GrindEQ__6_11_}, where $\theta \in \left\{\theta _{+} ,\theta _{-} \right\}$ and $\overline{\lambda}\in (0,1)$. Therefore, it can be said that the verification of the fulfillment of assumptions \eqref{GrindEQ__6_1_} and \eqref{GrindEQ__6_2_} in Theorem 6.1 is more constructive than the verification of the fulfillment of assumptions like \eqref{GrindEQ__6_11_}.

Applying the method of proof part (ii) of Theorem 5.1 and using Lemmas 3.2 and 4.2, it is easy to prove the following Theorem.

\textbf{Theorem 6.2.} Let the functions $L(\cdot )$, $L_{\dot{x}} (\cdot )$ and $L_{\dot{y}} (\cdot )$ be continuously differentiable with respect to their arguments and the function $\varphi \left(\cdot \right)-$ be twice continuously differentiable. In addition, let the admissible function is the weak local minimum of the problem \eqref{GrindEQ__1_1_}, \eqref{GrindEQ__1_2_}. Then there exists a number  $\delta >0$, at which the following statements are true:

(i) If the assumptions made in part (i) of theorem 6.1 are satisfied, then along the function $\bar{x}(\cdot )$ for each point $\left(\eta ,(\bar{\lambda }-1)^{-1} \bar{\lambda }\eta ,\, \bar{\lambda }\right)\in B_{\delta } (0)\times B_{\delta } (0)\times (0,1)$, satisfying condition \eqref{GrindEQ__6_1_} ((6.2)) the inequality \eqref{GrindEQ__6_3_} ((6.4)) holds;

(ii) If the assumptions made in part (ii) of theorem 6.1 are satisfied, then for each point $\left(\eta ,(\bar{\lambda }-1)^{-1} \bar{\lambda }\eta ,\, \bar{\lambda }\right)\in B_{\delta } (0)\times B_{\delta } (0)\times (0,1)$, satisfying condition \eqref{GrindEQ__6_5_}, the equality \eqref{GrindEQ__6_6_} holds.

\textbf{Remark 6.2.} If  $\theta\in \left[t_1-h,t_1\right),\theta\in \left(t_1-h,t_1\right]$ and  $\theta\in (t_1-h,t_1)$ the problem \eqref{GrindEQ__1_1_}, \eqref{GrindEQ__1_2_} is studied quite similarly to [23]. In addition, for the case $\theta\in [t_1-h,t_1]$, the statements of Theorems 5.1, 6.1 and 6.2 simplify in form, since by assumption the equality   $L(t,x,y,\dot{x},\dot{y})=0$,$ (t,x,y,\dot{x},\dot{y})\in (t_1,+\infty )\times R^n\times R^n\times R^n\times R^n$.

\section{Example and discussions}

To demonstrate the effectiveness of the obtained results, for example, the efficiency of Theorem 5.1, let's consider the following example.

 \textbf{Example 7.1}. Let's consider the problem:
\begin{equation} \label{GrindEQ__7_1_} 
S(x(\cdot ))=\int^3_0{\left[\left(1-x\right){\dot{x}}^2-\left(1+y\right){\dot{y}}^2+\dot{x}\dot{y}\right]dt}\to\mathop{min}_{x(\cdot )},  
\end{equation} 
\begin{equation} \label{GrindEQ__7_2_} 
x\left(t\right)=0,\ \ \ t\in \left[-1,0\right],\ \ x(3)=0, 
\end{equation} 
where $y=y\left(t\right)=x\left(t-1\right),h=1,L\left(t,x,y,\dot{x},\dot{y}\right)=\left(1-x\right){\dot{x}}^2-\left(1+y\right){\dot{y}}^2+\dot{x}\dot{y}$.

Let us examine the minimum admissible function $\overline{x}\left(t\right)=0,t\in \left[-1,3\right].$ Along this function, considering \eqref{GrindEQ__1_3_}, \eqref{GrindEQ__1_4_}, and \eqref{GrindEQ__1_6_}, we present the following calculations of the form:
\[\bar{L}\left(t\right)={\bar{L}}_x\left(t\right)={\bar{L}}_y\left(t\right)={\bar{L}}_{\dot{x}}\left(t\right)={\bar{L}}_{\dot{y}}\left(t\right)=0,t\in \left[0,3\right],\] 
\[{\bar{L}}_{\dot{y}}\left(t+1\right)={\bar{L}}_y\left(t+1\right)= 0, t\in \left[0,2\right];\] 
\[
\bar{L}\left( t,\xi ;\dot{\overline{x}}(\cdot )\right) =\xi ^{2},t\in \left[
0,3\right] ,\bar{L}\left( t+1,\xi ;\dot{\overline{y}}(\cdot )\right) =-\xi
^{2},\ \ t\in \left[ 0,2\right] ;
\]

\[
\bar{L}_{x}\left( t,\xi ;\dot{\overline{x}}(\cdot )\right) ={-\xi }^{2},\
t\in \left[ 0,3\right] ,\bar{L}_{y}\left( t+1,\xi ;\dot{\overline{y}}(\cdot
)\right) =-\xi ^{2},\ t\in \left[ 0,2\right] ;
\]%
\[
E\left( \bar{L}\right) \left( t,\xi ;\dot{\overline{x}}\left( \cdot \right)
\right) =\xi ^{2},t\in \left[ 0,3\right] ,\ \ E\left( \bar{L}\right) \left(
t+1,\xi ;\dot{\overline{y}}\left( \cdot \right) \right) =-\xi ^{2},\ \ t\in %
\left[ 0,2\right] ;
\]%
\[
M\left( \bar{L}_{x}\right) \left( t,\lambda ,\xi ;\dot{\overline{x}}\left(
\cdot \right) \right) =-\lambda \xi ^{2}-(1-\lambda )\left(\frac{\lambda }{
\lambda -1}\xi \right)^{2}=-\frac{\lambda }{1-\lambda }\xi ^{2}\ ,\ \ t\in \left[
0,3\right] ;
\]%
\[
M\left( \bar{L}_{y}\right) \left( t+1,\lambda ,\xi ;\dot{\bar{y}}\left(
\cdot \right) \right) =-\lambda \xi ^{2}-(1-\lambda )\left(\frac{\lambda }{%
\lambda -1}\xi \right)^{2}=-\frac{\lambda }{1-\lambda }\xi ^{2}\ ,\ \ t\in \left[
0,2\right] .
\]
    Based on these calculations, we come to the following conclusions:

(a) by virtue of \eqref{GrindEQ__1_7_} or \eqref{GrindEQ__1_8_}, the admissible function $\bar{x}\left(\cdot \right)=0$   is the extremal of the problem \eqref{GrindEQ__7_1_}, \eqref{GrindEQ__7_2_};

(b) the function $\bar{x}\left(\cdot \right)=0$ satisfies Weierstrass condition \eqref{GrindEQ__1_9_}:   \[\xi^2\ge 0,\forall \xi\in R,\forall t\in \left[2,3\right]\;\mbox{and} \;\xi^2-\xi^2=0,\forall t\in \left[0,2\right],\forall \xi\in R;\]

(c) it is evident that along the function $\bar{x}(\cdot )=0\ $ for all points  $(\lambda,\eta,\frac{\lambda}{\lambda-1}\eta)\in (0,1)\times R\times R$ of the Weierstrass condition \eqref{GrindEQ__1_9_} degenerates at any point in the interval (0,2).

Considering these  conclusions, we apply Theorem 5.1. As a result, statement \eqref{GrindEQ__5_2_} does not hold:
\[\frac{2\lambda}{\lambda-1}\eta^2=0\ ,\forall (t,\lambda,\eta)\in (0,2)\times (0,1)\times R\backslash \left\{0\right\}.\] 

 Therefore, by virtue of part (i) of Theorem 5.1, the extremal $\bar{x}(\cdot )=0\ $  cannot be a strong local minimum in the problem \eqref{GrindEQ__7_1_}, \eqref{GrindEQ__7_2_}.

Further, since equation \eqref{GrindEQ__5_2_} does not hold for all $\eta \neq 0$ and  $\lambda\in \left(0,1\right)\ $it is clear that by virtue of part (ii) of Theorem 5.1, the extremal $\bar{x}(\cdot )=0\ $ is not even a weak local minimum in the problem \eqref{GrindEQ__7_1_}, \eqref{GrindEQ__7_2_}. 

As far as we know, analogs of Theorems 5.1, 6.1 and 6.2 have not been obtained in the theory of singular optimal control.

Comparing this work with the work [23], we easily conclude:

(a) if in the problem \eqref{GrindEQ__1_1_}, \eqref{GrindEQ__1_2_} the number $h=0$, then the statements of Theorems 5.1, 6.1 and 6.2 coincide with the analogous results in [23];

(b) In this case, the approach proposed by us allows us to obtain analogues and other results of the work [23].

 (c) problem \eqref{GrindEQ__1_1_}, \eqref{GrindEQ__1_2_} is studied in a manner entirely analogous to [23], if Weierstrass conditions \eqref{GrindEQ__1_9_} degenerate at points of the segments $[t_1-h,t_1]$;

In conclusion, we also note that the application of the method outlined in this work is \linebreak promising for study more complex variational problems with delay such as variational problem with moving ends and a variational problem on a conditional extremum.

\end{document}